\pdfoutput=1 
\documentclass{amsart}
\usepackage[english]{babel}
\usepackage{amsfonts}
\usepackage{amsthm}
\usepackage{amsmath}
\usepackage{amssymb}
\usepackage{mathrsfs} 
\usepackage{enumerate}
\usepackage{mathtools}
\usepackage{tikz}
\usepackage{tikz-cd}
\usepackage[all]{xy}
\usepackage{graphicx}
\usepackage{amsmath}
\usepackage{adjustbox}
\usepackage{graphicx}
\usepackage{extpfeil}
\usepackage{comment}
\usepackage{mathabx}
\usepackage{float}
\usepackage[labelformat=empty]{caption}
\usepackage[toc]{appendix}
\usepackage{hyperref}
\usepackage{resizegather}
\usepackage[activate={true, nocompatibility}, final, tracking=true, kerning=true, spacing=true, factor=1100, stretch=10, shrink=10]{microtype}
\usepackage{dsfont}
\usetikzlibrary{matrix}
\usepackage{cleveref}
\usepackage{ulem}

\theoremstyle{plain}
\newtheorem{thm}{Theorem}[section]

\newtheorem{defn}[thm]{Definition}
\newtheorem{lemma}[thm]{Lemma}

\newtheorem{prop}[thm]{Proposition}

\newtheorem{remark}[thm]{Remark}

\newtheorem*{thm*}{Main Theorem}

\tikzset{
  symbol/.style={
    draw=none,
    every to/.append style={
      edge node={node [sloped, allow upside down, auto=false]{$#1$}}}
  }
}

\microtypecontext{spacing=nonfrench}

\newcommand{\on}{\operatorname}

\newcommand{\Fr}{\on{Fr}}

\newcommand{\cA}{\mathcal{A}}

\newcommand{\cC}{\mathcal{C}}

\newcommand{\cE}{\mathcal{E}}
\newcommand{\cF}{\mathcal{F}}
\newcommand{\cG}{\mathcal{G}}
\newcommand{\cH}{\mathcal{H}}
\newcommand{\cI}{\mathcal{I}}

\newcommand{\cL}{\mathcal{L}}
\newcommand{\cM}{\mathcal{M}}
\newcommand{\cN}{\mathcal{N}}
\newcommand{\cO}{\mathcal{O}}
\newcommand{\cP}{\mathcal{P}}

\newcommand{\cT}{\mathcal{T}}

\newcommand{\sE}{\mathscr{E}}
\newcommand{\sF}{\mathscr{F}}

\newcommand{\Ft}{\mathfrak{t}}

\newcommand{\bA}{\mathbb{A}}

\newcommand{\bF}{\mathbb{F}}

\newcommand{\bN}{\mathbb{N}}

\newcommand{\bR}{\mathbb{R}}
\newcommand{\bZ}{\mathbb{Z}}

\newcommand{\assa}
{\eva\mapsto\evalpha}

\newcommand{\assb}
{\evb\mapsto\evbeta}

\newcommand{\asstau}
{\evtau\mapsto\evbeta}

\newcommand{\good}{\mathrm{vg}}

\usepackage{color} 

\newcommand{\hig}{\mathcal{H}iggs}
\newcommand{\con}{\mathcal{C}onn}

\newcommand{\End}{\mathcal{E}nd}
\newcommand{\Bun}{\mathcal{B}un}
\newcommand{\Isom}{\mathcal{I}som}
\newcommand{\Aut}{\mathcal{A}ut}

\usepackage[toc]{appendix}

\newcommand{\eva}{\boldsymbol{a}}
\newcommand{\evb}{\boldsymbol{b}}

\newcommand{\evm}{\boldsymbol{m}}
\newcommand{\evalpha}{\boldsymbol{\alpha}}
\newcommand{\evbeta}{\boldsymbol{\beta}}
\newcommand{\evtau}{\boldsymbol{\tau}}
\newcommand{\evtheta}{\boldsymbol{\theta}}

\newcommand{\alphaop}{\alpha_{\operatorname{op}}^{\widetilde{\eva}}}
\newcommand{\betaop}{\beta_{\operatorname{op}}^{\widetilde{\evb}}}

\hypersetup{pdfborder=1 1 0.1} 

\title[Non-abelian Hodge for semistable parabolic connections in positive characteristic]{Non-abelian Hodge correspondence for semistable parabolic connections on curves in positive characteristic}
\author[Shiyu Shen]{Shiyu Shen}
\address{Institute of Science and Technology Austria (ISTA),
Klosterneuburg, Austria}
\email{shiyu.shen.math@gmail.com}

\begin{document}
\pagestyle{headings}
\begin{abstract}
    Let $C$ be a smooth projective curve over an algebraically closed field $k$ of characteristic $p>0$. We establish a non-abelian Hodge correspondence between semistable parabolic connections on $C$ and semistable parabolic Higgs bundles on its Frobenius twist $C'$, allowing arbitrary prescribed residue eigenvalues, parabolic weights and filtration types. The relationship between the parabolic weights and residue eigenvalues on the two sides is encoded in a table analogous to Simpson's table in \cite{simpson1990harmonic}. We obtain this result by extending the notion of ``very good regular connections'' to incorporate a residue-dependent degree condition, and by developing a Cartier descent formalism for parabolic connections. Working with filtered sheaves, we obtain a particularly clean description of this Cartier descent formalism that potentially extends to higher dimensional varieties.  
\end{abstract}

\maketitle

\section{Introduction}

\subsection{Non-abelian Hodge correspondence for Riemann surfaces}
Let $C$ be a compact Riemann surface. The non-abelian Hodge correspondence, established through the work of Hitchin \cite{hitchin1987self}, Donaldson \cite{donaldson1987twisted}, Corlette \cite{corlette1988flat}, and Simpson \cite{simpson1992higgs}, gives a homeomorphism between the moduli space of semistable Higgs bundles of degree zero and the moduli space of semistable flat connections. This groundbreaking work has led to major developments in related fields of mathematics, most notably the study of higher Teichm\"uller theory through Higgs bundles \cite{hitchin1992lie}, and the $P=W$ conjecture formulated in \cite{dCHM2012}, now proved for $\on{GL}_n$ in \cite{MaulikShenPW} and \cite{HMMSPW}. 

A logarithmic version of this non-abelian Hodge correspondence was established by Simpson \cite{simpson1990harmonic}. In order to incorporate arbitrary residue eigenvalues for logarithmic connections, one needs to introduce parabolic structures on the Higgs bundle side of the correspondence. Simpson's theorem produces a correspondence between parabolic connections and parabolic Higgs bundles with prescribed residue
eigenvalues, parabolic weights and filtration types. The relationship between the parabolic
weights and residue eigenvalues on the two sides is encoded in a table, see \cite[p.~720]{simpson1990harmonic}.

\subsection{Non-abelian Hodge correspondence for curves in positive characteristic}
Let $k$ be an algebraically closed field of characteristic $p>0$. Compared to the case over the complex numbers, a substantial feature of
the characteristic $p>0$ setting is that the algebra of crystalline differential operators has a large center. In fact, it is an Azumaya algebra (i.e. a sheaf of algebras \'etale-locally isomorphic to a matrix algebra) over the Frobenius twist of the cotangent bundle of the underlying space, see \cite{BMR08}. Based on this Azumaya property, Ogus and Vologodsky \cite{ov07} established a non-abelian Hodge correspondence assuming certain nilpotence and lifting assumptions. Over a smooth projective curve, Groechenig extended Tsen's theorem to show that the Azumaya algebra of differential operators splits on all spectral curves, thereby establishing a non-abelian Hodge correspondence over the entire Hitchin base; see \cite{groechenig-moduli-flat-connections}. The correspondence was generalized to arbitrary reductive groups in \cite{chen-zhu} under suitable characteristic assumptions. In \cite{dCGZ}, the notion of very good splittings
is introduced to obtain a correspondence between
semistable Higgs bundles and connections. The matching of semistable objects for general reductive groups is established in \cite{herrero-zhang-semistable}. 

The generalization of this non-abelian Hodge correspondence over curves to the logarithmic setting was initiated in \cite{shen2018tamely}.  Working over $\on{GL}_n$, \cite{shen2018tamely} established a correspondence under a certain regularity condition on the residue eigenvalues for logarithmic connections, see \Cref{defn: regular tuple}. With the same regularity assumption, the matching of semistable objects is established in \cite{de2025logarithmic}, using a similar ``very good'' notion as in \cite{dCGZ}. Similar to the characteristic zero setting, going beyond this regularity assumption requires the introduction of parabolic structures on the Higgs bundle side, and corresponding results appear in \cite{li-sun2021tame} and \cite{barz2025logarithmic}. 

\subsection{Main results}
Let $C$ be a smooth projective curve over an algebraically closed field $k$ of characteristic $p>0$. All theorems in this paper admit straightforward generalizations to any reduced divisor $D\subset C$. We restrict to the case of a single point $q\in C$ for simplicity of exposition. The main goal of this paper is to establish a non-abelian Hodge correspondence between rank $n$ semistable parabolic connections on $C$ and rank $n$ semistable
parabolic Higgs bundles on its Frobenius twist $C'$, allowing arbitrary prescribed residue
eigenvalues, parabolic weights and filtration types, with no restriction on $p$. The relationship between the parabolic
weights and residue eigenvalues on the two sides is encoded in the following table analogous to Simpson's
table in \cite{simpson1990harmonic}. 
\begin{center}
\begin{tabular}{|c|c|c|}
\hline
 & Parabolic Higgs bundle on $C'$
 & Parabolic connection on $C$ \\
\hline
Parabolic weight
 & $\alpha=(\beta-\nu(b))/p$ & $\beta$ \\
\hline
Residue eigenvalue
 & $a=b^p-b$ & $b$ \\
\hline
\end{tabular}
\end{center}
Here $\nu$ is a fixed map from $k$ to $\{-p+1, \dots, -1,0\}$ whose restriction to each fiber of the Artin--Schreier map $b\mapsto b^p-b$ is a bijection; see \Cref{sec: very good} for its precise definition. Since parabolic weights lie in $[0,1)$, the above table offers a bijection between pairs $(a, \alpha)$ and $(b, \beta)$ in $k\times [0,1)$. 

Before giving the precise statement of our main theorem, we introduce the following notation. Let $\eva$ be an unordered $n$-tuple of elements in $k$, and let $\evalpha$ be a set of parabolic weights. We introduce the notion of an assignment $\assa$, which assigns a parabolic weight in $\evalpha$ to each entry of $\eva$, counted with multiplicity; see \Cref{def: assignment}. Let $\hig^{\assa}(C)$ (resp. $\con^{\assa}(C)$) be the moduli stack of parabolic Higgs bundles (resp. connections) whose residue at $q$
acts on each graded quotient with eigenvalues
prescribed by $\assa$; see \Cref{def: substack according to assignment}. Let $\assa$ and $\assb$ be two assignments related through the table above; see \Cref{sec: Statement of the main theorem} for the precise relation. The moduli stack $\hig^{\assa}(C')$ admits the Hitchin morphism $h'$, and $\con^{\assb}(C)$ admits the de Rham--Hitchin morphism $h_{\on{dR}}$ to the same Hitchin base:
\[
\xymatrix{
\con^{\assb}(C)\ar[dr]^{h_{\on{dR}}}&& \hig^{\assa}(C')\ar[dl]_{h'}\\
&A_{\eva}.
}
\]
Our main theorem establishes an isomorphism between $\con^{\assb}(C)$ and $\hig^{\assa}(C')$ up to twisting by a certain Picard torsor, and  parabolic semistability is preserved under this correspondence. 

Let $\evb$ be a regular tuple in the sense of \Cref{defn: regular tuple} that maps to $\eva$ under the componentwise Artin--Schreier map. The moduli stack $\cH_{\evb}$ of regular\footnote{By regular we mean that the $p$-curvature of the connection is regular as a twisted Higgs field.} logarithmic connections with eigenvalues of residue $\evb$ introduced in \cite{shen2018tamely} plays a crucial role in the logarithmic non-abelian Hodge correspondence per se, as well as in subsequent works \cite{li-sun2021tame},  \cite{barz2025logarithmic} and \cite{de2025logarithmic}. In particular, it was shown in \cite{shen2018tamely} that when $p>n$, $\cH_{\evb}$ maps surjectively to the corresponding Hitchin base $A_{\eva}$, and is a torsor for the relative Picard stack $\cP_{\eva}$ of spectral curves over $A_{\eva}$. We provide a slightly different local construction in this paper to remove this restriction on the characteristic; see \Cref{prop: regular connection is torsor}. 

In order to match semistable objects in the parabolic setting, we generalize the notion of ``very good regular connections'' introduced in \cite{dCGZ} and \cite{de2025logarithmic} by incorporating a residue-dependent degree condition; see \Cref{defn: good component}. We show that this ``very good'' substack $\cH_{\evb}^{\mathrm{vg}}\subset\cH_{\evb}$ maps surjectively to $A_{\eva}$ and is a $\cP_{\eva}^o$-torsor, where $\cP_{\eva}^o\subset \cP_{\eva}$ is the substack parametrizing line bundles having degree zero on each integral component; see \Cref{prop: very good is torsor}. Up to canonical equivalence via Hecke transformations, this $\cP_{\eva}^o$-torsor $\cH_{\evb}^{\mathrm{vg}}$ is independent of the choice of a regular Artin--Schreier preimage $\evb$, and we denote it by $\cH_{\eva}^{\mathrm{vg}}$. 

Now we are ready to state the main theorem of the paper. Note that the tuple $\evb$ underlying $\assb$ need not be regular. Let $\hig^{\assa,\mathrm{ss}}(C',d)\subset \hig^{\assa}(C')$ (resp. $\con^{\assb,\mathrm{ss}}(C,d)\subset\con^{\assb}(C)$) be the substack classifying $\evalpha$-semistable (resp. $\evbeta$-semistable) objects of parabolic degree $d$. The following is proved in \Cref{thm: stacky main para} and \Cref{thm: main para}.
\begin{thm}\label{thm: main thm intro}
    Let $\assa$ and $\assb$ be two assignments related through the table above. There exists an isomorphism of $A_{\eva}$-stacks  
\[ \cH_{\eva}^{\good}\times_{A_{\eva}}^{\cP_{\eva}^o} \hig^{\assa}(C')
         \cong \con^{\assb}(C). 
\]
For any $d\in\bR$, this isomorphism restricts to an isomorphism
\[  \cH_{\eva}^{\good}\times_{A_{\eva}}^{\cP_{\eva}^o} \hig^{\assa, \mathrm{ss}}(C',d)
         \cong \con^{\assb,\mathrm{ss}}(C, pd).
\] 
\end{thm}

A key ingredient in the proof of \Cref{thm: main thm intro} is the establishment of a Cartier descent formalism for parabolic connections, building on work of Li and Sun~\cite{li-sun2021tame}. Let $\evtau$ be an unordered $n$-tuple of elements in $\bF_p$. We denote by $\con_0^{\asstau}(C)\subset \con^{\asstau}(C)$ the substack of connections with zero $p$-curvature, and by $\Bun^{\evalpha}(C', \evm')$ the moduli stack of parabolic bundles with parabolic weights $\evalpha$ and filtration type $\evm'$. The following is proved in \Cref{prop: logarithmic Cartier descent}. 

\begin{prop}\label{prop: logarithmic Cartier descent intro}
For any assignment $\asstau$, there exists an isomorphism
\[\con_0^{\asstau}(C)\cong \Bun^{\evalpha}(C', \evm'), 
\]  
where the parabolic weights $\evalpha$ and filtration type $\evm'$ arise from the assignment $\assa$ constructed from $\asstau$ through the table above. 
\end{prop}

A crucial step towards the proof of \Cref{prop: logarithmic Cartier descent intro} is to show that every object in $\con_0^{\asstau}(C)$ is locally isomorphic to a certain model object; see \Cref{lemma: local isomorphism}. This is a parabolic analogue of
\cite[Lemma 5.21]{li-sun2021tame}, and we offer an independent proof using Lagrange interpolation polynomials. 

Utilizing the correspondence between parabolic bundles and filtered sheaves, we obtain a particularly clean description of
this Cartier descent formalism, as stated in the following theorem proved in \Cref{thm: filtered sheaf cartier descent}. 
\begin{thm}\label{thm: filtered sheaf cartier descent intro}
Under the isomorphism in \Cref{prop: logarithmic Cartier descent intro}, a filtered connection $\sE=(\{\sE_{x}\}_{x\in \bR}, \nabla_{\sE})$ corresponds to the filtered sheaf $\sF=\{\sF_{x}\}_{x\in \bR}$ defined by
\begin{equation}\label{eqn: filteres sheaf descent intro}
\sF_{x}\coloneqq(\on{Fr}_C)_*\ker\bigl(
\nabla_{\sE}\colon \sE_{px}
\longrightarrow
\sE_{px}\otimes\omega_C(q)\bigr),
\end{equation}
where $\on{Fr}_C\colon C\to C'$ is the relative Frobenius.
\end{thm}

The inverse of this construction is provided in \Cref{Remark: inverse of Phi}. Note that this formulation of Cartier descent for filtered sheaves through (\ref{eqn: filteres sheaf descent intro}) treats all choices of residue eigenvalues in $\bF_p$ and parabolic weights simultaneously, and the description involves neither root stacks nor parahoric group schemes. This formulation also potentially extends
to higher dimensional varieties. We plan to investigate this generalization in future work.

\subsection{Organization of the paper}
In \Cref{sec: Logarithmic Higgs bundles and logarithmic connections}, we review various constructions related to logarithmic Higgs bundles and connections. We show the substack of regular connections maps surjectively to the corresponding Hitchin base, and is a torsor for the relative Picard stack of spectral curves.  In \Cref{sec: very good}, we first introduce the residue-dependent degree condition that defines the substack $\cH_{\evb}^{\mathrm{vg}}$ of ``very good logarithmic connections''. Then we show $\cH_{\evb}^{\mathrm{vg}}$ is a $\cP^o_{\eva}$-torsor independent of the choice of the Artin--Schreier preimage $\evb$, up to canonical equivalence. In \Cref{sec: Parabolic structures and residue assignments}, we first review the definition of parabolic bundles, Higgs bundles, and connections, as well as their parabolic semistability conditions. Then we introduce the notion of an assignment of parabolic weights to eigenvalues. In \Cref{sec: main theorem}, we state and prove the main parabolic non-abelian Hodge correspondence \Cref{thm: main thm intro} of the paper, after introducing the analogue of Simpson's table that encodes the relationship between the parabolic weights and
residue eigenvalues on the two sides. We also establish the parabolic Cartier descent formalism \Cref{prop: logarithmic Cartier descent intro}. In \Cref{sec: filtered sheaf Cartier descent}, working with filtered sheaves, we provide another description of this Cartier descent formalism by proving \Cref{thm: filtered sheaf cartier descent intro}. 

\subsection{Acknowledgments}
The author thanks Siqing Zhang for contributing to the proof of \Cref{prop: very good is torsor}, and for helpful discussions on related subjects. The author was partially supported by the European Union's Horizon 2020 research and innovation program under the Marie Skłodowska-Curie grant agreement No. 101034413 and the Austrian Science Fund (FWF), grant 10.55776/P35847.

\subsection{Declaration of AI usage}
The author used GPT-5.6 sol (very high), GPT-6 (very high) and Claude Opus 5.5 (extra) to proofread and improve the exposition of the manuscript. The manuscript was written entirely by the author, who takes full responsibility for its content.

\section{Logarithmic Higgs bundles and logarithmic connections}\label{sec: Logarithmic Higgs bundles and logarithmic connections}

\subsection{Logarithmic Higgs bundles}
Let $k$ be an algebraically closed field. Let $C$ be a smooth projective curve over $k$, and let $q\in C$ be a marked point. Let $\omega_C$ be the sheaf of differentials on $C$. A logarithmic Higgs bundle over $(C, q)$ is a pair $(E,\phi)$ where $E$ is a vector bundle over $C$ and $\phi$ is an $\cO_C$-linear morphism 
\[\phi\colon E\longrightarrow E\otimes\omega_C(q).
\]
Let $\hig(C,n)$ be the moduli stack of rank $n$ logarithmic Higgs bundles over $(C,q)$. For a $k$-scheme $S$, $\hig(C,n)(S)$ is the groupoid of pairs $(E,\phi)$, where $E$ is a rank $n$ vector bundle on $C\times S$ and $\phi$ is an $\cO_{C\times S}$-linear morphism
\[
\phi\colon E\longrightarrow E\otimes p_C^*\omega_C(q),
\]
where $p_C\colon C\times S\to C$ is the projection map. By taking the coefficients of the characteristic polynomial of $\phi$, we obtain the so-called Hitchin morphism, which we denote by 
\[h\colon \hig(C,n)\longrightarrow A(\omega_C(q))\coloneqq\bigoplus_{i=1}^n H^0(C,\omega_C^i(iq)).
\]
\subsubsection{Residues of logarithmic Higgs bundles}
Recall that the fiber $\omega_C(q)|_q$ carries a canonical trivialization induced by the residue map\[\omega_C(q)|_q\xrightarrow[]{\sim}\cO_q\cong k.\]
For any choice of a local coordinate $x$ of $C$ vanishing at $q$, this trivialization maps the section $\frac{dx}{x}$ to $1$.
Let $(E,\phi)$ be a logarithmic Higgs bundle over $(C,q)$. The residue of $\phi$ at $q$
\[\on{res}_q(\phi)\in H^0(q, \End(E)\otimes \omega_C(q))\cong H^0(q,\End(E|_q))
\] 
is the restriction of $\phi$ to $q$, where the second isomorphism is induced by the residue trivialization.

Let $\mathfrak{t}\subset \mathfrak{gl}_n$ be the Cartan subalgebra given by diagonal matrices, and let $W$ be the Weyl group acting on $\mathfrak{t}$. Consider the evaluation-at-$q$ morphism 
\begin{equation}\label{eqn: evalution of Hitchin base at q}
    \on{ev}_q\colon A(\omega_C(q))=\bigoplus_{i=1}^n H^0(C,\omega_C^i(iq))\longrightarrow \bigoplus_{i=1}^n H^0(q,\omega_C^i(iq))\xrightarrow[]{\sim}\Gamma(q,\cO_q^{\oplus n})\cong \mathfrak{t}(k)/W,
\end{equation}
where the isomorphism $\bigoplus_{i=1}^n H^0(q,\omega_C^i(iq))\xrightarrow[]{\sim}\Gamma(q,\cO_q^{\oplus n})$ is induced by the tensor powers of the residue trivialization $\omega_C(q)|_q\xrightarrow[]{\sim}\cO_q$, and the last isomorphism is given by the standard identification $\mathbb{A}^n\cong\mathfrak{t}/W$ via the coefficients of the characteristic polynomial. The Hitchin morphism induces a morphism \[\Gamma(q, \End(E|_q))\longrightarrow \Gamma(q,\cO_q^{\oplus n})\cong \mathfrak{t}(k)/W
\] 
given by sending an endomorphism to the coefficients of its characteristic polynomial, which determine and are determined by the eigenvalues of the endomorphism. We parametrize the closed points of $\mathfrak{t}/W$  by unordered $n$-tuples of elements in $k$, written as 
\begin{equation}\label{eqn: veceva}
\eva=a_1^{n_1}a_2^{n_2}\cdots a_s^{n_s},
\end{equation}
with $a_j\in k, n_j\in \bN_{>0}$ and $\sum_{j=1}^{s}n_j=n$.

\subsection{Logarithmic connections}
Now we assume $k$ is of characteristic $p>0$. A logarithmic connection over $(C,q)$ is a pair $(E,\nabla)$ where $E$ is a vector bundle over $C$ and $\nabla: E\to E\otimes_{\cO_C}\omega_C(q)$ is a $k$-linear map that satisfies the Leibniz rule.
Let $\con(C,n)$ be the moduli stack of rank $n$ logarithmic connections over $(C,q)$. For a $k$-scheme $S$, $\con(C,n)(S)$ is the groupoid of pairs $(E,\nabla)$, where $E$ is a rank $n$ vector bundle on $C\times S$ and $\nabla$ is an $\cO_{S}$-linear morphism
\[
\nabla\colon E\longrightarrow E\otimes p_C^*\omega_C(q)
\]
that satisfies the Leibniz rule. 

\subsubsection{The de Rham--Hitchin morphism}
Let $\cT_C$ be the tangent sheaf of $C$. Let $C'=C\times_{k,F_k}k$ be the Frobenius twist of $C$, and let $\Fr_C\colon C\to C'$ be the relative Frobenius morphism. Let $q'=\Fr_C(q)\in C'$. For any logarithmic connection $(E, \nabla)$, we associate with $\nabla$ its $p$-curvature
\[\psi_{\nabla} :\on{Fr}_C^{*}(\mathcal{T}_{C'}(-q'))\longrightarrow \mathcal{E}nd(E),
\]
which is the $\mathcal{O}_C$-linear morphism corresponding to the
$p$-linear morphism
\[
\psi_\nabla:\mathcal T_C(-q)\longrightarrow\mathcal End(E),
\qquad
\psi_\nabla(\partial)
=
\nabla_\partial^p-\nabla_{\partial^{[p]}}
\]
for any $\partial\in\mathcal{T}_C(-q)(U)$ and $U\subseteq C$ open. 
The $p$-curvature $\psi_{\nabla}$ can also be viewed as a twisted Higgs field
\[ \psi_{\nabla}: E\longrightarrow E\otimes \on{Fr}_C^{*}(\omega_{C'}(q')).
\]
The coefficients of the characteristic polynomial of $\psi_\nabla$ define a point $a\in A\bigl(\Fr_C^*(\omega_{C'}(q'))\bigr)$. By the same argument as in \cite[Proposition~3.2]{laszlo-pauly-p-curvature}, this point lies in the image of the pullback map
$\on{Fr}_C^*:
A\bigl(\omega_{C'}(q')\bigr)\to A\bigl(\Fr_C^*(\omega_{C'}(q'))\bigr)$. Thus we obtain a morphism
 \[h_{\on{dR}}: \con(C,n)\longrightarrow A(\omega_{C'}(q'))\]
 which we call the de Rham--Hitchin morphism.

\subsubsection{Residues of logarithmic connections}\label{sec: res conn}
Given a logarithmic connection $(E,\nabla)$ over $(C, q)$, the residue trivialization of $\omega_C(q)$ gives rise to the residue $\on{res}_q(\nabla)\in\Gamma(q, \End(E|_q))$, and the Frobenius pullback of the residue trivialization of $\omega_{C'}(q')$ gives rise to the residue $\on{res}_q(\psi_{\nabla})\in\Gamma(q, \End(E|_q))$. 
We have the identity of residues in $\Gamma(q,\End(E|_q))$
\[\on{res}_q(\psi_\nabla)=\on{res}_q(\nabla)^p-\on{res}_q(\nabla),\]
see \cite[Lemma~4.1]{shen2018tamely}.
Thus we have the identity of eigenvalues in $\mathfrak{t}(k)/W$:
\begin{equation}\label{eqn: as on eval}   \eva(\on{res}_q(\psi_\nabla))=\eva(\on{res}_q(\nabla))^p-\eva(\on{res}_q(\nabla)).
\end{equation}
Let $\on{AS}: \mathfrak{t}(k)/W\to\mathfrak{t}(k)/W $ be the morphism given by $\eva\mapsto \eva^p-\eva$, which is induced by the \'etale componentwise Artin--Schreier morphism on $\bA^{n}$ and the quotient morphism by the natural action of the permutation group $W=S_n$.

\subsection{Spectral description of logarithmic connections}\label{sec: torsor conn}
For the rest of the paper, we write $A\coloneqq A(\omega_{C'}(q'))$. Let $S$ be the universal spectral curve over $A$.
Let $\cP(S/A)$ be the relative Picard stack for $S/A$.
Let $\cT F(S/A)$ be the $A$-stack parametrizing torsion free rank 1 sheaves on $S$ over $A$ in the sense of Schaub \cite{schaub-courbes-spectrales}, see also \cite[\S2.4]{de2017support-sln}.
By the Beauville--Narasimhan--Ramanan (BNR) correspondence as in \cite[Propositions 2.1 and 5.1]{schaub-courbes-spectrales}, there is a canonical isomorphism of $A$-stacks \[\hig(C',n)\cong \cT F(S/A),\] mapping a logarithmic Higgs bundle to its spectral sheaf. 

Let $S^{\Fr}
\coloneqq
C\times_{\Fr_C,C'}S$ be the Frobenius pullback of the universal spectral curve $S$. Let $\cP(S^{\Fr}/A)$ be the relative Picard stack for $S^{\Fr}/A$. Let $\cT F(S^{\Fr}/A)$ be the $A$-stack parametrizing torsion free rank 1 sheaves on $S^{\Fr}$ over $A$. Taking the spectral sheaf of the $p$-curvature $\psi_\nabla$, we obtain a canonical morphism 
of $A$-stacks
\begin{equation}\label{eqn: drbnr sheaf}
    \con(C,n)\to \cT F(S^{\Fr}/A).
\end{equation}
Let \[\cH\coloneqq\con(C,n)\times_{\cT F(S^{\Fr}/A)} \cP(S^{\Fr}/A)\] be the open substack of $\con(C,n)$ characterized by the spectral sheaf of $\psi_{\nabla}$ being invertible. To characterize $\cH$, we introduce the following:
\begin{defn}\label{defn: regular tuple}
    Let $\evb=b_1^{n_1}b_2^{n_2}\cdots b_s^{n_s}\in \mathfrak{t}(k)/W$ be parametrized as in (\ref{eqn: veceva}). We say $\evb$ is a regular tuple if for each pair $j,j'$ with $j\neq j'$, we have $b_j-b_{j'}\notin \bF_p$.
\end{defn}
Given a connection $(E,\nabla)$ with eigenvalues of residue $\evb$,
if $\evb$ is not regular, then by writing residues in Jordan normal forms, the $p$-curvature $\psi_\nabla$ is not regular in the Lie-theoretic sense, i.e. the corresponding Higgs bundle $(E,\psi_\nabla)$ does not land in the open Picard stack inside the Higgs moduli, see \cite[Proposition 2.6.1.(1)]{chzh17} and \cite[\S2.2.8]{de2025logarithmic}. 
By \cite[Proof of Lemma 3.15]{de2025logarithmic}, objects in $\cH$ can only have regular eigenvalues of residues. For any regular tuple $\evb$, we define $\cH_{\evb}\subset \cH$ to be the substack characterized by the condition that the residue of the connection has eigenvalues $\evb$. Let $\eva=\on{AS}(\evb)$, and let $A_{\eva}\coloneqq \on{ev}_{q'}^{-1}(\eva)\subset A$, where $\on{ev}_{q'}$ is the evaluation-at-$q'$ morphism (\ref{eqn: evalution of Hitchin base at q}). It follows from (\ref{eqn: as on eval}) that $h_{\on{dR}}(\cH_{\evb})\subseteq A_{\eva}$. 

For the rest of the paper, we write $\cP\coloneqq \cP(S/A)$ and $\cP_{\eva}\coloneqq \cP(S/A_{\eva})$. Recall that by \cite[Lemma 4.11]{shen2018tamely}, each $\cH_{\evb}$ admits a $\cP_{\eva}$-action, and becomes a pseudo $\cP_{\eva}$-torsor under this action, see \cite[Proposition 4.13]{shen2018tamely}. We now show the following
\begin{prop}\label{prop: regular connection is torsor}
    For any regular tuple $\evb$, $\cH_{\evb}$ is a $\cP_{\eva}$-torsor. 
\end{prop}
\begin{proof}
   It is enough to show that $\cH_{\evb}\to A_{\eva}$ is surjective. When $\on{char}(k)=p$ satisfies $p>n$, this is shown in \cite[Proposition 4.6]{shen2018tamely}. The construction of local models in Step 2 of the proof of \cite[Proposition 4.6]{shen2018tamely} relies on the assumption $p>n$. We now provide a slightly different argument that does not rely on this assumption. Let $\evb=b_1^{n_1}b_2^{n_2}\cdots b_s^{n_s}$ be parametrized as in (\ref{eqn: veceva}), and let $\on{AS}(\evb)=\eva=a_1^{n_1}a_2^{n_2}\cdots a_s^{n_s}$, with $\on{AS}(b_j)=a_j$. Let $D_q=\on{Spec}(k[[x]])$ be the formal disc around $q\in C$, and set $y=x^p\partial_x^p$. Let $a\in A_{\eva}$. When restricted to $D_q$, the Frobenius pullback of the spectral curve decomposes as 
\[S_a^{\Fr}\times_C D_q=\coprod_{j=1}^s S_{a,a_j}^{\Fr}. 
\]
Each $S_{a,a_j}^{\Fr}$ is of the form 
\[
S_{a,a_j}^{\Fr}\cong \on{Spec}\bigl(k[[x]][y]/(f_j)\bigr),
\]
where $f_j\in k[[x^p]][y]$ and satisfies $\bar{f}_j=(y-a_j)^{n_j}\in k[[x^p]][y]/x^pk[[x^p]][y]$. Our goal is to construct for each $j$ a logarithmic connection $(E_j, \nabla_j)$ over $D_q$ that satisfies the following three conditions:
\begin{enumerate}
    \item \label{surj: first} the eigenvalues of $\on{res}_q(\nabla_j)$ are all $b_j$,
	\item \label{surj: second}the spectral curve of $\psi_{\nabla_j}$ is $S_{a,a_j}^{\Fr}$,
	\item \label{surj: third}the spectral sheaf of $\psi_{\nabla_j}$ is an invertible sheaf on $S_{a,a_j}^{\Fr}$. 
\end{enumerate}
Since $S_{a,a_j}^{\Fr}$ is the spectrum of a local ring, (\ref{surj: third}) implies that the spectral sheaf of $\psi_{\nabla_j}$ is isomorphic to $\mathcal{O}_{S_{a,a_j}^{\Fr}}$. Let $e$ be its generator. A logarithmic connection with spectral curve $S_{a,a_j}^{\Fr}$ is determined by the connection acting on $e$, which can be written as $\nabla(e)=ge \frac{dx}{x}$, $g\in k[[x]][y]/(f_j)$. Assuming $g\in k[[x^p]][y]/(f_j)$,
conditions (\ref{surj: first}) and (\ref{surj: second}) are equivalent to 
\begin{equation}\label{eqn: local connection}
g^p-g=y,\qquad
g\equiv b_j\pmod{(x^p,y-a_j)}.
\end{equation}
Since $b_j$ is a simple root of $z^p-z=a_j$, the equations (\ref{eqn: local connection}) admit a unique solution by Hensel's lemma. Indeed, this solution can be written explicitly as 
\[
g=b_j-\sum_{m=0}^{\infty}(y-a_j)^{p^m}.
\]
The rest of the argument proceeds as in Steps 1 and 3
of the proof of \cite[Proposition 4.6]{shen2018tamely},
using Beauville--Laszlo gluing to obtain the desired
global logarithmic connection.
\end{proof}

\section{The stack of very good logarithmic connections}\label{sec: very good}
Let $\evb\in \mathfrak{t}(k)/W$ be a regular tuple as defined in \Cref{defn: regular tuple}. The goal of this section is to define the substack of so-called very good logarithmic connections $\cH_{\evb}^{\good}\subset \cH_{\evb}$, which plays a crucial role in our main theorem on the non-abelian Hodge correspondence for semistable parabolic connections.

Before defining the very good substack, we fix the following notation and conventions. Let $\sigma:k\to\mathbb F_p$ be an $\mathbb F_p$-linear retraction  of the inclusion $\mathbb F_p\hookrightarrow k$. We fix this choice throughout this paper.
Let $\tau: \bF_p\to \bZ$ be the set-theoretic section of the quotient map $\bZ\to \bF_p$ with image in $\{-p+1,-p+2,\dots,-1,0\}$.
Let \[\nu\coloneqq\tau\circ \sigma: k\to \bZ\] be the composition. Since $\sigma$ is an $\bF_p$-linear retraction, $\nu(b+c)=\tau(\sigma(b)+c)$ for $c\in\bF_p$, so $\nu$ restricts to a bijection from each Artin--Schreier fiber $b+\bF_p$ onto $\{-p+1,\dots,-1,0\}$.

Let $\evb=b_1^{n_1}b_2^{n_2}\cdots b_s^{n_s}$ be parametrized as in (\ref{eqn: veceva}), and let $\on{AS}(\evb)=\eva=a_1^{n_1}a_2^{n_2}\cdots a_s^{n_s}$, with $\on{AS}(b_j)=a_j$. It follows from \Cref{prop: regular connection is torsor} that $h_{\on{dR}}(\cH_{{\evb}})=A_{\eva}$. Let $a\in A_{\eva}$. Given an object $\cE\coloneqq(E,\nabla)$ in the fiber groupoid $(\cH_{\evb})_a$, let $L_{\cE}$ be the spectral sheaf of the $p$-curvature $\psi_{\nabla}$, i.e. the line bundle on $S^{\Fr}_{a}$ given by the image of the morphism \eqref{eqn: drbnr sheaf}.
The corresponding spectral curve $S^{\Fr}_{a}$ decomposes into integral components 
\[S^{\Fr}_{a}=\sum_{t=1}^N n_t S^{\Fr}_{a,t},\]
where $S^{\Fr}_{a,t}$ are the distinct reduced irreducible components and $n_t$ are their multiplicities.
Let $\pi^{\Fr}: S^{\Fr}_{a}\to C$ be the projection.
We denote
\[(\pi^{\Fr})^{-1}(q)\times_{S^{\Fr}_{a}} S^{\Fr}_{a,t} =\sum_{l=1}^{M_{t}} c_{t,l} a_{t,l}\] 
with $c_{t,l}\in \bN_{>0}$ and $a_{t,l}$ a closed point in the fiber $S^{\Fr}_{a,t}\times_C q\subset\on{Tot}\bigl(\on{Fr}_C^*(\omega_{C'}(q'))\bigr)\times_C q$. The latter
is canonically identified with $\mathbb{A}_k^1$ through the Frobenius pullback of the residue trivialization, so that we can write $a_{t,l}\in k$ unambiguously.
Since $\evb$ is a regular tuple, there is a unique $b_j$ in the decomposition $\evb=b_1^{n_1}b_2^{n_2}\cdots b_s^{n_s}$ satisfying \[b_j^p-b_j=a_{t,l}\in k.\]
We write this $b_j$ as $b_{t,l}$.

\begin{defn}\label{defn: good component}
We say that $\cE$ is very good if, for each integral component
$S^{\Fr}_{a,t}$ of $S^{\Fr}_{a}$, the following equality holds in $\mathbb{Z}$:
\begin{equation}\label{eqn: the degrees}  
    \deg_{S^{\Fr}_{a,t}}   \bigl(L_{\cE}|_{S^{\Fr}_{a,t}}\bigr)
    +\sum_{l=1}^{M_t} c_{t,l}\nu(b_{t,l})=0. 
\end{equation}
\end{defn}

\begin{defn}
 Let $\cH_{\evb}^{\good}$ be the substack of $\cH_{\evb}$ classifying logarithmic connections $\mathcal{E}$ that restrict to a very good logarithmic connection in the sense of \Cref{defn: good component} at each geometric point of the test scheme.   
\end{defn}
Let $\cP^o_{\eva}\subseteq\cP_{\eva}$ be the substack characterized by
the condition that, at each geometric point of the test
scheme with image $a\in A_{\eva}$, the corresponding line bundle on $S_a$ has degree zero
on each integral component of the spectral curve $S_a$.
\begin{prop}\label{prop: very good is torsor}
    The $\cP_{\eva}$-action on $\cH_{\evb}$ restricts to a $\cP^o_{\eva}$-action on $\cH^{\good}_{\evb}$. Moreover, $\cH^{\good}_{\evb}$ is a $\cP^o_{\eva}$-torsor. 
\end{prop}
\begin{proof}
Let $a\in A_{\eva}$. Let $W:S_a^{\Fr}=C\times_{\Fr_C,C'}S_a\to S_a$ and $\pi^{\Fr}:S_a^{\Fr}=C\times_{\Fr_C,C'}S_a\to C$
be the projections. Let $\cE=(E,\nabla_E)\in \cH_{\evb}$ such that $h_{\on{dR}}(\cE)=a\in A_{\eva}$, and let $L_{\cE}$ be the spectral sheaf of its $p$-curvature. Let $L'$ be a line bundle on $S_a$. The $\cP_{\eva}$-action on $\cH_{\evb}$ maps $(\cE, L')$ to
\[
\bigl(\pi^{\Fr}_*(L_{\cE}\otimes W^*L'), \nabla_{\otimes}\bigr),
\]
where $\nabla_{\otimes}$ is the tensor product connection of $\nabla_E$ acting on $L_{\cE}$ and the canonical connection on $W^*L'$, see \cite[Lemma 4.9]{shen2018tamely}. Now we assume $\cE\in \cH_{\evb}^{\good}$. Note that the logarithmic connection resulting from the action of $L'$ on $\cE$ remains in $\cH_{\evb}^{\good}$ if and only if $W^*L'$ has degree zero on each integral component $S_{a,t}^{\Fr}$ of $S_a^{\Fr}$. Let $S_{a,t}^{\mathrm{red}}$ be the integral
component of $S_a$ corresponding to $S_{a,t}^{\Fr}$
under the universal homeomorphism
$W:S_a^{\Fr}\to S_a$. Then we have $\bigl(C\times_{\Fr_C,C'}S_{a,t}^{\mathrm{red}}\bigr)^{\mathrm{red}}
=S_{a,t}^{\Fr}$, therefore
\[
\on{deg}(W^*L'|_{S_{a,t}^{\Fr}})=0\Leftrightarrow \on{deg}(L'|_{S^{\mathrm{red}}_{a,t}})=0. 
\]
Since all integral components of $S_a$ arise this way for some $t$, it follows that $\cH_{\evb}^{\good}$ is a pseudo $\cP^o_{\eva}$-torsor. 

Now we show that $\cH_{\evb}^{\good}$ maps surjectively to $A_{\eva}$. Fix $a\in A_{\eva}$. By \Cref{prop: regular connection is torsor}, there exists $\cE=(E,\nabla_E)\in \cH_{\evb}$ such that $h_{\on{dR}}(\cE)=a$. Denote $Z_t\coloneqq C\times_{\Fr_C,C'}S_{a,t}^{\mathrm{red}}\subseteq S^{\Fr}_{a}$. Let $\pi_t: S_{a,t}^{\mathrm{red}}\to C'$ and $\pi^{\Fr}_t: Z_t\to C$ be the projection maps. Note that $Z_t=e_tS^{\Fr}_{a,t}$ for some $e_t\in \bN_{>0}$. Since $Z_t$ is defined through Frobenius pullback, its defining ideal is preserved by the connection $\nabla_E$. Therefore $\nabla_E$ restricts to a logarithmic connection on $E_t\coloneqq (\pi^{\Fr}_t)_*(L_{\cE}|_{Z_t})$ whose residue has eigenvalues given by $\prod_{l=1}^{M_t}b_{t,l}^{e_tc_{t,l}}$. By the residue--degree formula, we have
\[
\on{deg}(E_t)+\sum_{l=1}^{M_t}e_tc_{t,l}b_{t,l}=0 \quad\text{in }k. 
\]
It follows that $\sum_{l=1}^{M_t}e_{t}c_{t,l}b_{t,l}$ lies in $\mathbb F_p$.
Since $\sigma$ is an $\mathbb F_p$-linear retraction and
$\nu=\tau\circ\sigma$, we have
\[
\sum_{l=1}^{M_t}e_tc_{t,l}\nu(b_{t,l})
\equiv \tau(\sum_{l=1}^{M_t}e_tc_{t,l}b_{t,l})\pmod p.
\]
By Riemann--Roch, we have 
\[
\on{deg}(E_t)=\on{deg}(L_{\cE}|_{Z_t})+p\on{deg}((\pi_t)_*\cO_{S_{a,t}^{\mathrm{red}}}).
\]
It follows that 
\[
\on{deg}(L_{\cE}|_{Z_t})+ \sum_{l=1}^{M_t}e_{t}c_{t,l}\nu(b_{t,l})\equiv 0 \pmod p.
\]
Set $d_t=
\frac{1}{p}\left(
\on{deg}(L_{\cE}|_{Z_t})
+\sum_{l=1}^{M_t}e_tc_{t,l}\nu(b_{t,l})
\right)$, and let $L'$ be a line bundle on $S_a$ with degree $-d_t$ on each integral component $S_{a,t}^{\mathrm{red}}$. Let $\cE'$ be the logarithmic connection arising from the action of $L'$ on $\cE$. Since $L_{\cE'}=L_{\cE}\otimes W^*L'$, we have
\[
\on{deg}(L_{\cE'}|_{Z_t})+ \sum_{l=1}^{M_t}e_{t}c_{t,l}\nu(b_{t,l})=0 \quad\text{in }\bZ. 
\]
Combining with 
\[
\on{deg}(L_{\cE'}|_{Z_t})=e_t \on{deg}(L_{\cE'}|_{S_{a,t}^{\Fr}}),
\]
we obtain the desired equality (\ref{eqn: the degrees}) in \Cref{defn: good component}. 
\end{proof}

\begin{prop}\label{prop: heck transform for good connections}
    Let $\evb'\in \Gamma(q,\cO_q^{\oplus n})\cong \mathfrak{t}(k)/W$ be another regular tuple such that $\on{AS}(\evb')=\on{AS}(\evb)=\eva$. 
    Then we have an isomorphism  $\mathcal{H}ecke_{\evb', \evb}:\cH_{\evb}^{\good}\xrightarrow{\sim}\cH_{\evb'}^{\good}$ of $\cP_{\eva}^o$-torsors induced by Hecke transformation. Given a third regular tuple $\evb''$ mapping to $\eva$ under $\on{AS}$, we have a canonical $2$-isomorphism
    \[\mathcal{H}ecke_{\evb'',\evb'}\circ\mathcal{H}ecke_{\evb',\evb}\simeq\mathcal{H}ecke_{\evb'',\evb}.\]
\end{prop}

\begin{proof}
Let $\evb=b_1^{n_1}b_2^{n_2}\cdots b_s^{n_s}$, and let $\on{AS}(\evb)=\eva=a_1^{n_1}a_2^{n_2}\cdots a_s^{n_s}$, with $\on{AS}(b_j)=a_j$. After possibly reordering the indices, $\evb'$ can be written as $\evb'=(b_1')^{n_1}(b_2')^{n_2}\cdots (b_s')^{n_s}$ with $\on{AS}(b_j')=\on{AS}(b_j)=a_j$. Let $\cE=(E,\nabla)\in \cH_{\evb}$ such that $h_{\on{dR}}(\cE)=a\in A_{\eva}$. The spectral correspondence for $(E, \psi_{\nabla})$ gives rise to a line bundle $L_{\mathcal{E}}$ on $S_a^{\Fr}$ such that $\pi^{\Fr}_{*}L_{\mathcal{E}}\cong E$. Let $x$ be a local parameter of $C$ at $q$, so that $\widehat{\cO}_{C,q}\cong k[[x]]$ and $K_q\cong k((x))$. Let $D_q\coloneqq \on{Spec} \widehat{\cO}_{C,q}\cong \on{Spec}k[[x]]$ be the formal disc around $q\in C$, and let $D^{\times}_q\coloneqq \on{Spec} K_q\cong \on{Spec}k((x))$ be the punctured disc. By definition of $A_{\eva}$, we have
\[(\pi^{\Fr})^{-1}(q)=\sum_{j=1}^s n_ja_j.
\]
It follows that $S_a^{\Fr}\times_C D_q$ decomposes as a disjoint union
\[S_a^{\Fr}\times_C D_q=\coprod_{j=1}^s S_{a,a_j}^{\Fr}.
\]
Therefore $E|_{D_q}$ decomposes as
\[
E|_{D_q}=\bigoplus_{j=1}^s M_j,
\]
with $M_j$ supported on $S_{a,a_j}^{\Fr}$ for each $j$. Note that the connection $\nabla$ preserves each $M_j$, and $\on{res}_q(\nabla)$ has eigenvalue $b_j$ when acting on $M_j$. Consider the $k[[x]]$-lattice $M'\subset E|_{D_q^{\times}}$ defined by
\begin{equation}\label{eqn: hecke modifications}
M'= x^{\nu(b_1')- \nu(b_1)}M_1\oplus\cdots \oplus x^{\nu(b_s')- \nu(b_s)}M_s.
\end{equation}
By the Beauville--Laszlo theorem \cite{BeauvilleLaszlo1995}, $E|_{C\setminus q}$ and $M'$ glue to a new vector bundle $E'$ on $C$ through their canonical identification over $D_q^{\times}$. The logarithmic connection $\nabla$ on $E$ induces a logarithmic connection on $E'$, which we still denote by $\nabla$. Now we
define 
\[
\mathcal{H}ecke_{\evb', \evb}: \cH_{\evb}\longrightarrow \cH_{\evb'}, \qquad (E, \nabla)\mapsto (E', \nabla).
\]
Note that for each $1\leq j\leq s$, $\on{res}_q(\nabla)$ has eigenvalue $\nu(b_j')-\nu(b_j)+b_j=b_j'$ when acting on $x^{\nu(b_j')- \nu(b_j)}M_j$.  When restricted to $D_q\subset C$, each integral component $S_{a,t}^{\Fr}$ of $S_{a}^{\Fr}$ decomposes as
\[
S_{a,t}^{\Fr}\times_{C}D_q=\coprod_{l=1}^{M_t} S_{a,t, a_{t,l}}^{\Fr}.
\]
Let $L_{\mathcal{E},t,l}$ denote the restriction of $L_{\mathcal{E}}$ to $S_{a,t, a_{t,l}}^{\Fr}$,  then the modification \eqref{eqn: hecke modifications} corresponds to
\[
L_{\mathcal{E},t,l}\mapsto L_{\mathcal{E},t,l}\bigl((\nu(b_{t,l})- \nu(b_{t,l}'))\on{div}_{S_{a,t, a_{t,l}}^{\Fr}}(x)\bigr)=L_{\mathcal{E},t,l}\bigl((\nu(b_{t,l})- \nu(b_{t,l}'))c_{t,l}a_{t,l}\bigr). 
\]
It follows that $\mathcal{H}ecke_{\evb', \evb}$ maps $\cH_{\evb}^{\good}$ to $\cH_{\evb'}^{\good}$, and is equivariant with respect to the $\cP^{o}_{\eva}$-action. An inverse of $\mathcal{H}ecke_{\evb', \evb}$ is given by $\mathcal{H}ecke_{\evb, \evb'}$, therefore we have an isomorphism
\[
\mathcal{H}ecke_{\evb', \evb}:\cH_{\evb}^{\good}\xrightarrow{\sim}\cH_{\evb'}^{\good} 
\]
of $\cP^{o}_{\eva}$-torsors. The composition $2$-isomorphism follows directly from the definition.
\end{proof}

It follows that, up to canonical isomorphism, the $\cP_{\eva}^{o}$-torsor $\cH_{\evb}^{\good}$ is independent of the choice of the regular Artin--Schreier preimage $\evb$. From now on we denote this torsor by $\cH_{\eva}^{\good}$. 

\section{Parabolic structures and residue assignments}\label{sec: Parabolic structures and residue assignments}
\subsection{Parabolic Higgs bundles and connections}\label{sec: parabolic notation}

Let $\evm=(m_1, m_2, \dots, m_r)$ be an ordered tuple of positive integers that sum to $n$. A quasi-parabolic bundle of filtration type $\evm$ is a vector bundle $E$ together with a filtration of $k$-vector spaces 
\[0=E_{0}\subset E_{1}\subset \cdots\subset E_{r-1}\subset E_{r}=E_{q}\coloneqq E|_q,\]
such that $\dim_k(E_{j}/E_{j-1})=m_j$ for $j=1, 2,\dots,r$. We denote such a filtration by $E_{\bullet}$.

A quasi-parabolic Higgs bundle is a triple $(E, E_{\bullet}, \phi)$, where $(E,E_{\bullet})$ is a quasi-parabolic bundle and $\phi$ is a logarithmic Higgs field whose residue $\on{res}_q(\phi)$ preserves the filtration $E_{\bullet}$: $\on{res}_{q}(\phi)(E_{j})\subseteq E_{j}$. 
Similarly, a quasi-parabolic connection is a triple $(E, E_{\bullet}, \nabla)$, where $(E, E_{\bullet})$ is a quasi-parabolic bundle and $\nabla$ is a logarithmic connection whose residue $\on{res}_q(\nabla)$ preserves the filtration $E_{\bullet}$: $\on{res}_{q}(\nabla)(E_{j})\subseteq E_{j}$.
Let $\hig(C,\evm)$ (resp. $\con(C,\evm)$) denote the moduli stack of quasi-parabolic Higgs bundles (resp. connections) of filtration type $\evm$.

A parabolic bundle $(E, E_{\bullet}, \evalpha)$ is a quasi-parabolic bundle $(E, E_{\bullet})$ together with a set of parabolic weights $\evalpha\subset [0,1)$, enumerated in decreasing order as
\[
\evalpha=(\alpha_1,\alpha_2,\ldots,\alpha_r), \qquad
1>\alpha_1>\alpha_2>\cdots>\alpha_r\geq 0.
\]
Let $\Bun^{\evalpha}(C, \evm)$, $\hig^{\evalpha}(C,\evm)$, and $\con^{\evalpha}(C,\evm)$ denote the moduli stacks of parabolic bundles, Higgs bundles, and connections with parabolic weights $\evalpha$ and filtration type $\evm$, respectively.

\begin{remark}\label{remark: parabolic bundle and parabolic group torsor}
Let $\evtheta\in \mathfrak{t}(\bR)$. For a root $\gamma$ of $\on{GL}_n$, set
\[
m_\gamma \coloneqq
\left\lceil-\gamma(\evtheta)\right\rceil=
\min\bigl\{m\in\mathbb Z \mid \gamma(\evtheta)+m\geq 0\bigr\}. 
\]
We denote by $\cG_{\evtheta}$ the smooth affine group scheme over $C$
which is isomorphic to $\on{GL}_n$ over $C\setminus q$ and whose group of $k[[x]]$-points is generated by $T(k[[x]])$ together with the root subgroups
\[
U_\gamma\bigl(x^{m_\gamma}k[[x]]\bigr).
\]
Equivalently,
\[
\cG_{\evtheta}(R[[x]])
=
\left\{
(g_{ij})\in\on{GL}_n(R((x)))
\;\middle|\;
\begin{array}{l}
g_{ij}\in x^{\lceil\theta_j-\theta_i\rceil}R[[x]]
\quad\text{for all }i,j,\\
\det(g)\in R[[x]]^\times
\end{array}
\right\}.
\]
Different orderings of $\evtheta$ give isomorphic group schemes, with the isomorphism induced by conjugation by the corresponding permutation matrix. Hence the isomorphism class of $\cG_{\evtheta}$ only depends on the class of $\evtheta$ in $\mathfrak{t}(\bR)/W$.

Now let
\[
\evalpha=(\alpha_1,\alpha_2,\dots,\alpha_r),
\qquad
\evm=(m_1,m_2,\dots,m_r),
\]
with $1>\alpha_1>\alpha_2>\cdots>\alpha_r\geq 0$ and
$\sum_{j=1}^r m_j=n$. We associate to $(\evalpha,\evm)$ the element
\[
\evtheta({\evalpha,\evm})
=
\on{diag}\bigl(
\alpha_1^{(m_1)},
\alpha_2^{(m_2)},
\dots,
\alpha_r^{(m_r)}
\bigr)
\in \mathfrak t(\mathbb{R}),
\]
where $\alpha_j^{(m_j)}$ means that $\alpha_j$ is repeated $m_j$ times. Let $\evtheta=(\theta_1,\theta_2,\dots,\theta_n)\in \mathfrak{t}(\mathbb{R})$ be an element in the same Weyl group orbit of $\evtheta({\evalpha,\evm})$. A parabolic vector bundle on $C$ with parabolic weights
$\evalpha$ and filtration type
$\evm$ is equivalent to a
$\cG_{\evtheta}$-torsor on $C$. We describe this correspondence as follows. Fix a model object $\cF^0=(F^0, F^0_{\bullet})$, where $F^0=\bigoplus_{i=1}^n \cO_C$, and the parabolic structure is obtained by assigning the parabolic weight $\theta_i$ to the $i$-th summand $(\cO_C)_q\subset F^0_q$. The sheaf of automorphisms ${\cA}ut(\cF^0)$ is equal to $\cG_{\evtheta}$. For any parabolic bundle $\cF=(F, F_{\bullet})$ with parabolic weights $\evalpha$ and filtration type $\evm$, we associate with it the sheaf $\cI som(\cF^0, \cF)$ of isomorphisms from $\cF^0$ to $\cF$. Note that $\cI som(\cF^0, \cF)$ is an ${\cA}ut(\cF^0)$-torsor through precomposition, and the evaluation map induces
\[
{\cI}som(\cF^0,\cF)\times^{{\cA}ut(\cF^0)}\cF^0\cong  \cF.
\]
\end{remark}

For each $1\leq j\leq r$, we denote by
$\on{gr}_{\alpha_j}(E_q)$ the graded piece of weight $\alpha_j$, i.e. $\on{gr}_{\alpha_j}(E_q)
\coloneqq
E_j/E_{j-1}$. The $\evalpha$-parabolic degree of $(E, E_{\bullet})$ is given by 
\[\on{pardeg}(E,E_{\bullet},\evalpha)\coloneqq\deg(E)+\sum_{j=1}^{r} \alpha_{j}m_j.\]
Given a subbundle $E'$ of $E$, we obtain a parabolic subbundle $(E', E_{\bullet}', \evalpha')$ of $(E,E_{\bullet}, \evalpha)$ as follows. The filtration $E_{\bullet}'$ is given by intersecting $E_q'$ with the filtration $E_{\bullet}$. The induced parabolic weights $\evalpha'$ on $(E', E_{\bullet}')$ are given by $\alpha'_{k}=\max\{\alpha_{j}\mid E'_{k}\subseteq E_{j}\}$. We can then make sense of the $\evalpha$-parabolic slope (semi)stability of parabolic bundles, Higgs bundles, and connections, where in the latter two cases we consider subbundles preserved by the Higgs field or connection, respectively. Let $\Bun^{\evalpha\text{-}\mathrm{ss}}(C,\evm)$, $\hig^{\evalpha\text{-}\mathrm{ss}}(C,\evm)$, and  $\con^{\evalpha\text{-}\mathrm{ss}}(C,\evm)$ denote the moduli stacks of $\evalpha$-semistable parabolic bundles, Higgs bundles, and connections of filtration type $\evm$, respectively.

\subsection{Residue assignments and the associated moduli substacks}

Let $\evalpha$ be a set of parabolic weights on a quasi-parabolic bundle of filtration type $\evm$ as in \Cref{sec: parabolic notation}. In order to obtain a non-abelian Hodge correspondence for parabolic Higgs bundles $(E,E_{\bullet}, \evalpha, \phi)$, we introduce notation and conventions to fix the multiset of eigenvalues for the action of $\on{res}_q(\phi)$ on each graded piece of the filtration $E_{\bullet}$. We fix $\eva\in \mathfrak{t}(k)/W$, and fix an (arbitrary) lift of $\eva$ to $\widetilde{\eva}\in \mathfrak{t}(k)$, which is represented by an ordered tuple $\widetilde{\eva}=(a_1, a_2, \dots, a_n)$. For each $1\leq i\leq n$, we assign a parabolic weight $\alpha_{\operatorname{op}}^{\widetilde{\eva}}(i)\in \evalpha$ to each $a_{i}$ in $\widetilde{\eva}$, subject to the requirement that each $\alpha_j\in \evalpha$ appears exactly $m_j$ times. We denote such an assignment of parabolic weights to eigenvalues by $\eva\mapsto \evalpha$. More precisely, 
\begin{defn}\label{def: assignment}
    Given a filtration type $\evm$, a set of parabolic weights $\evalpha\subset [0,1)$
with size equal to the length of $\evm$, and an element $\eva\in \Ft(k)/W$, an assignment $\assa$ (of parabolic weights to eigenvalues) of filtration type $\evm$ is defined as an equivalence class of maps
\[
\alpha_{\operatorname{op}}^{\widetilde \eva}\colon
\{1,\dots,n\}\longrightarrow \evalpha,
\qquad
i\longmapsto \alpha_{\operatorname{op}}^{\widetilde \eva}(i),
\]
subject to
\[
\#\bigl\{i\mid\alpha_{\operatorname{op}}^{\widetilde \eva}(i)=\alpha_j\bigr\}=m_j
\]
for each $j$. Two such maps are equivalent if they differ by a
permutation of indices preserving $\widetilde{\eva}$.
Whenever we write $\alpha_{\operatorname{op}}^{\widetilde{\eva}}$,
we choose a representative of the equivalence class. 
\end{defn} 

Now we define the following substacks of $\hig^{\evalpha}(C,\evm)$ and $\con^{\evalpha}(C,\evm)$ characterized by the condition that the multiset of eigenvalues of the residue acting on each graded piece is given by $\assa$. These substacks are the main players in the non-abelian Hodge correspondences for parabolic Higgs bundles established in this paper. 

\begin{defn}\label{def: substack according to assignment}
Let $\eva\mapsto \evalpha$ be an assignment of filtration type $\evm$ as in \Cref{def: assignment}.  

\begin{enumerate}
    \item We denote by $\hig^{\assa}(C)\subseteq\hig^{\evalpha}(C,\evm)$ the substack classifying parabolic Higgs bundles $(E,E_{\bullet},\phi)$ that satisfy the following extra condition on the residue $\on{res}_q(\phi)$ of the Higgs field. For each $1\leq j \leq r$, we require that the multiset of eigenvalues of
$\on{res}_q(\phi)$ acting on $\on{gr}_{\alpha_j}(E_q)$ is given by
\[
\bigl\{a_i \mid \alpha_{\operatorname{op}}^{\widetilde{\eva}}(i)=\alpha_j\bigr\},
\]
counted with multiplicity. Denote $\hig^{\assa,\mathrm{ss}}(C)\coloneqq \hig^{\eva\mapsto\evalpha}(C)\cap\hig^{\evalpha\text{-}\mathrm{ss}}(C,\evm)$.
\item 
We denote by $\con^{\eva\mapsto \evalpha}(C)\subseteq\con^{\evalpha}(C,\evm)$ the substack classifying parabolic connections $(E,E_{\bullet},\nabla)$ that satisfy the following extra condition on the residue $\on{res}_q(\nabla)$ of the connection. For each $1\leq j \leq r$, we require that the multiset of eigenvalues of
$\on{res}_q(\nabla)$ acting on $\on{gr}_{\alpha_j}(E_q)$ is given by
\[
\bigl\{a_i \mid \alpha_{\operatorname{op}}^{\widetilde{\eva}}(i)=\alpha_j\bigr\},
\]
counted with multiplicity. Denote $\con^{\assa,\mathrm{ss}}(C)\coloneqq \con^{\eva\mapsto\evalpha}(C)\cap\con^{\evalpha\text{-}\mathrm{ss}}(C,\evm)$.
\end{enumerate}
\end{defn}

\section{Non-abelian Hodge correspondence for semistable parabolic connections}\label{sec: main theorem}

\subsection{Statement of the main theorem}\label{sec: Statement of the main theorem}

Let $\eva\in \Ft(k)/W$ and let $\evalpha$ be a set of parabolic weights. Fix an assignment $\eva\mapsto\evalpha$ as in \Cref{def: assignment}. 
Since the restriction of $\nu$ to each fiber of the Artin--Schreier map is a bijection onto $\{-p+1,\dots,-1,0\}$, the following lemma is elementary:
\begin{lemma}
There is a unique ordered tuple $\widetilde{\evb}=(b_1, b_2,\dots, b_n)\in k^n$ such that 
\begin{enumerate}
    \item $b_{i}^p-b_{i}=a_{i}$.
    \item $0\leq p\alphaop(i)+\nu(b_{i})<1$. 
\end{enumerate}
\end{lemma}
The real numbers $p\alphaop(i)+\nu(b_{i})$ make up a new set of parabolic weights $\evbeta$. We obtain an assignment $\assb$ from $\assa$ by setting
\begin{equation}\label{eqn: assb from assa} 
\betaop(i)=p\alphaop(i)+\nu(b_{i}).   
\end{equation}
Conversely, given an assignment $\assb$, we can obtain $\assa$ by setting
\begin{equation} \label{eqn: assa from assb}    \widetilde{\eva}=\on{AS}(\widetilde{\evb})  \qquad\text{and}\qquad
\alphaop(i)=\frac{\betaop(i)-\nu(b_i)}{p}.
\end{equation}
Note that (\ref{eqn: assb from assa}) and (\ref{eqn: assa from assb}) are inverse to each other, thus $\assa$ and $\assb$ uniquely determine each other.

For any $d\in \bR$, let $\hig^{\assa, \mathrm{ss}}(C,d)\subset\hig^{\assa, \mathrm{ss}}(C)$ and $\con^{\assb, \mathrm{ss}}(C,d)\subset\con^{\assb, \mathrm{ss}}(C)$ be the substacks characterized by the underlying parabolic  bundle having parabolic degree $d$ at each geometric point of the test scheme. Now we are ready to state the main theorem of the paper.

\begin{thm}\label{thm: main para}
Let $\assa$ and $\assb$ be two assignments related through (\ref{eqn: assa from assb}). For any $d\in\bR$, there exists an isomorphism of $A_{\eva}$-stacks  
\[  \cH_{\eva}^{\good}\times_{A_{\eva}}^{\cP_{\eva}^o} \hig^{\assa, \mathrm{ss}}(C',d)
         \cong \con^{\assb,\mathrm{ss}}(C, pd).
\] 
\end{thm}
    If we choose the trivial filtration type $\evm=(n)$ and trivial parabolic weight $\evalpha=\{0\}$, then we have $\evbeta=\{0\}$ and $\evb$ is a regular tuple. In this case,  we recover the semistable logarithmic non-abelian Hodge correspondence as in \cite[Theorem 3.24]{de2025logarithmic} over $A_{\eva}$.

\subsection{Parabolic non-abelian Hodge without semistability}
Before proving \Cref{thm: main para}, we state and prove the following non-abelian Hodge correspondence for parabolic connections without imposing semistability. 
\begin{thm}\label{thm: stacky main para}
Let $\assa$ and $\assb$ be two assignments related through (\ref{eqn: assa from assb}). There exists an isomorphism  of $A_{\eva}$-stacks  
\[ \cH_{\eva}^{\good}\times_{A_{\eva}}^{\cP_{\eva}^o} \hig^{\assa}(C')
         \cong \con^{\assb}(C).
\]
\end{thm}

In \cite[Theorem~5.25]{li-sun2021tame}, Li and Sun constructed an isomorphism\footnote{\cite[Theorem~5.25]{li-sun2021tame} is proved assuming $p>n$. The general statement follows by combining their argument with \Cref{prop: regular connection is torsor}.}
\begin{equation}\label{eqn: Li-Sun naH}  \cH\times_A^{\mathcal P^{\text{ext}}}\mathfrak{X}^{\evbeta, \evm}\xrightarrow{\sim} \con^{\evbeta}(C, \evm).
\end{equation}
Here the stack $\mathfrak{X}^{\evbeta, \evm}$ classifies $(E, E_{\bullet},\nabla, \psi)$, where $(E, E_{\bullet},\nabla)$ is a parabolic connection with parabolic weights $\evbeta$ and filtration type $\evm$ such that its $p$-curvature is zero, and $\psi: E\to E\otimes \Fr_C^{*}\bigl(\omega_{C'}(q')\bigr)$ is an $\mathcal{O}_C$-linear morphism satisfying 
\begin{enumerate}
    \item\label{cond: horizontal} The morphism $\psi$ is horizontal with respect to the canonical connection on $\mathcal{H}om\bigl(E, E\otimes {\Fr}_C^{*}(\omega_{C'}(q'))\bigr)$ constructed from $\nabla$ on $E$ and the canonical connection on ${\Fr}_C^{*}(\omega_{C'}(q'))$.
    \item \label{cond: psi respects flag} The residue $\on{res}_q(\psi)$ respects the flag structure $E_{\bullet}$.
\end{enumerate}
The stack $\mathcal P^{\text{ext}}$ classifies $(E,\nabla, \psi)$, where $(E,\nabla)$ is a logarithmic connection with zero $p$-curvature, and $\psi: E\to E\otimes \Fr_C^{*}\bigl(\omega_{C'}(q')\bigr)$ is an $\mathcal{O}_C$-linear morphism that satisfies the horizontal condition (\ref{cond: horizontal}) above, as well as the extra requirement that the spectral sheaf of $\psi$ is invertible. The stacks $\mathfrak{X}^{\evbeta, \evm}$ and $\cP^{\text{ext}}$ map to the Hitchin base $A$ by taking the coefficients of the characteristic polynomial of $\psi$. 

The equivalence in \cite{li-sun2021tame} is described in terms of regular centralizers and parahoric group schemes. We give a more explicit description using spectral curves and parabolic structures as follows. Let $\cE=(E, \nabla_E)\in \cH$ such that $h_{\on{dR}}(\cE)=a\in A$, and let $L_{\cE}$ be the spectral sheaf of its $p$-curvature $\psi_{\nabla_E}$. Let $\cM=(M, M_{\bullet},\nabla_M, \psi)\in \mathfrak{X}^{\evbeta, \evm}$ that maps to the same $a\in A$, and let $\cF_{\psi}$ be the spectral sheaf of $\psi$. 
We define
\begin{equation}\label{eqn: stacky parabolic naH without torsor}
\cH\times_A\mathfrak{X}^{\evbeta, \evm}\longrightarrow\con^{\evbeta}(C, \evm)
\end{equation}
\[
(\cE, \cM)\mapsto \left(\pi_*^{\Fr}(L_\cE\otimes_{S_a^{\Fr}} \cF_{\psi}), \bigl(\pi_*^{\Fr}(L_\cE\otimes_{S_a^{\Fr}} \cF_{\psi})\bigr)_{\bullet}, \nabla_{\otimes}\right),
\]
where $\nabla_{\otimes}$ is the tensor product connection, and the parabolic structure $\bigl(\pi_*^{\Fr}(L_\cE\otimes_{S_a^{\Fr}} \cF_{\psi})\bigr)_{\bullet}$ is induced from the parabolic structure on $\pi_*^{\Fr}\cF_{\psi}$. 

Let $\cN=(N, N_{\bullet}, \nabla_N)\in \con^{\evbeta}(C, \evm)$ with $h_{\on{dR}}(\cN)=a$, and let $\cF_{\cN}$ be the spectral sheaf of its $p$-curvature $\psi_{\nabla_N}$. We define a morphism in the opposite direction
\begin{equation}\label{eqn: stacky parabolic naH without torsor: inverse}
  \cH\times_A\con^{\evbeta}(C, \evm)\longrightarrow \mathfrak{X}^{\evbeta, \evm} 
\end{equation}
\[
(\cE, \cN)\mapsto \left(\pi_*^{\Fr}(L^{-1}_\cE\otimes_{S_a^{\Fr}} \cF_{\cN}), \bigl(\pi_*^{\Fr}(L^{-1}_\cE\otimes_{S_a^{\Fr}} \cF_{\cN})\bigr)_{\bullet}, \nabla_{\otimes}, \psi\right).
\]
Here $\nabla_{\otimes}$ is the tensor product connection, with the connection on $L_{\cE}^{-1}$ defined by the natural connection on ${\cH}om$. The parabolic structure $\bigl(\pi_*^{\Fr}(L^{-1}_\cE\otimes_{S_a^{\Fr}} \cF_{\cN})\bigr)_{\bullet}$ is induced from the parabolic structure on $\pi_*^{\Fr}\cF_{\cN}$. The twisted Higgs field $\psi$ is given by the action of the tautological section on $S_a^{\Fr}$. Similar constructions to those in (\ref{eqn: stacky parabolic naH without torsor}) and (\ref{eqn: stacky parabolic naH without torsor: inverse}) endow $\mathcal{P}^{\on{ext}}$ with the structure of a commutative group stack. They also define actions of $\mathcal{P}^{\on{ext}}$ on $\mathfrak{X}^{\evbeta,\evm}$ and ${\cC}onn^{\evbeta}(C,\evm)$, with respect to which (\ref{eqn: stacky parabolic naH without torsor}) becomes $\mathcal{P}^{\on{ext}}$-equivariant, giving rise to the desired morphism (\ref{eqn: Li-Sun naH}). Combining \cite[Lemma~5.18]{li-sun2021tame} and \Cref{prop: regular connection is torsor}, we have that $\cH$ is a $\mathcal{P}^{\on{ext}}$-torsor, and it follows immediately that (\ref{eqn: Li-Sun naH}) is an isomorphism of stacks over $A$. 

Note that the $p$-curvature of $(E, \nabla)$ being zero implies that the residue $\on{res}_{q}(\nabla)$ is semisimple and has eigenvalues in $\mathbb{F}_p$. An element $\evtau\in \mathfrak{t}(\mathbb{F}_p)/W$ corresponds to an unordered $n$-tuple of elements in $\mathbb{F}_p$. We denote by $\mathfrak{X}^{\evbeta, \evm, \evtau}\subset \mathfrak{X}^{\evbeta, \evm}$ the substack classifying $(E, E_{\bullet},\nabla, \psi)$ with $\on{res}_q(\nabla)$ conjugate to $\evtau$. Then we have
\[
\mathfrak{X}^{\evbeta, \evm}=\coprod_{\evtau\in \mathfrak{t}(\mathbb{F}_p)/W} \mathfrak{X}^{\evbeta, \evm, \evtau}.
\]

Let $\asstau$ be an assignment of filtration type $\evm$. We define $\mathfrak{X}^{\asstau}$ to be the fiber product
\[
\mathfrak{X}^{\asstau}
\coloneqq
\mathfrak{X}^{\evbeta,\evm, \evtau}
\times_{\con^{\evbeta}(C,\evm)}
\con^{\asstau}(C),
\]
where the morphism
$\mathfrak{X}^{\evbeta,\evm, \evtau}
\to\con^{\evbeta}(C,\evm)$
is given by forgetting $\psi$. Then we have
\[
\mathfrak{X}^{\evbeta, \evm, \evtau}=\coprod\mathfrak{X}^{\asstau},
\]
where the disjoint union ranges over all assignments $\asstau$ of filtration type $\evm$.

\subsubsection{Correspondence between $\mathfrak{X}^{\asstau}$ and parabolic Higgs bundles}\label{sec: proof of identifying with parabolic higgs}

Now we describe $\mathfrak{X}^{\asstau}$. Let $\assa$ be the assignment constructed from $\asstau$ through (\ref{eqn: assa from assb})
\begin{equation}\label{eqn: asstau to assa}
\widetilde{\eva}=\on{AS}(\widetilde{\evtau})=\boldsymbol{0}  \qquad\text{and}\qquad
\alphaop(i)=\frac{\beta_{\on{op}}^{\widetilde{\evtau}}(i)-\nu(\tau_i)}{p}, 
\end{equation}
and we denote its filtration type by $\evm'$. A crucial intermediate result towards the proof of \Cref{thm: stacky main para} is the following
\begin{thm}\label{thm: identifying with parabolic higgs}
There exists an isomorphism 
$\mathfrak{X}^{\asstau}\cong \hig^{\evalpha}(C', \evm')$. This isomorphism is $\cP$-equivariant, where $\cP$
acts by tensoring spectral sheaves with line bundles
on the spectral curves, using Frobenius pullback
on the $\mathfrak{X}^{\asstau}$ side, with the parabolic structures transported accordingly.
\end{thm}

Let $\con_0^{\asstau}(C)\subset \con^{\asstau}(C)$ denote the substack classifying connections with vanishing $p$-curvature. Equivalently, $\con_0^{\asstau}(C)$ is obtained from $\mathfrak{X}^{\asstau}$ by forgetting the datum $\psi$. We fix a model object $\mathcal{E}^0=(E^0, E^0_{\bullet},\nabla^0)\in \con_0^{\asstau}(C)(k)$ defined as follows. Let $\widetilde{\evtau}=(\tau_1, \tau_2, \dots, \tau_n)\in\mathfrak{t}(\bF_p)$ be a lift of $\evtau\in\mathfrak{t}(\bF_p)/W$. Let $E^0=\bigoplus_{i=1}^n \mathcal{O}_C\bigl(-\nu(\tau_i)q\bigr)$, let $\nabla^0$ be the logarithmic connection on $E^0$ induced from the canonical connection on $\bigoplus_{i=1}^n \mathcal{O}_C$, and let $E^0_{\bullet}$ be the partial flag structure on $E^0_q$ obtained by assigning the parabolic weight $\beta_{\on{op}}^{\widetilde{\evtau}}(i)$ to the summand $\mathcal{O}_C(-\nu(\tau_i)q)_q\subset E^0_q$. Let $S$ be a $k$-scheme, and let $\mathcal{E}\in \con_0^{\asstau}(C)(S)$. We denote by $\cE_S^0$ the pullback of $\cE^0$ along the projection $p_C: S\times C\to C$. Let $\Isom(\cE_S^0, \cE)$ denote the sheaf of isomorphisms from $\cE_S^0$ to $\cE$ over $S\times C$. The sheaf $\Isom(\cE_S^0,\cE)$ carries a natural right action of the sheaf $\Aut(\cE_S^0)$ of automorphisms of $\cE_S^0$, given by
precomposition:
\[
\varphi\cdot g\coloneqq \varphi\circ g,
\qquad
\varphi\in\Isom(\cE_S^0,\cE),\quad
g\in\Aut(\cE_S^0).
\]

The following lemma is a parabolic analogue of
\cite[Lemma 5.21]{li-sun2021tame};
see also \cite[Lemma 5.26]{li-sun2021tame}.
We give a direct proof in our setting.
\begin{lemma}\label{lemma: local isomorphism}
Let $\cE\in \con_0^{\asstau}(C)(S)$. Then $\Isom(\cE_S^0,\cE)$ is an $\Aut(\cE_S^0)$-torsor.   
\end{lemma}
\begin{proof}
It is enough to show that there exists a Zariski open cover $\{U_i\to S\times C\}_{i\in I}$
such that $\cE|_{U_i}$ and $\cE_S^0|_{U_i}$ are isomorphic as objects of
$\con_0^{\asstau}(C)$ for every $i\in I$. Away from $S\times q$, the claim follows from Cartier descent: connections with vanishing $p$-curvature are obtained by Frobenius pullback of vector bundles on the Frobenius twist, which are Zariski locally trivial; see \cite[Theorem~5.1]{Katz1970}. Now we construct a Zariski open neighborhood $V\subset S\times C$ of $S\times q$, together with an isomorphism 
\[
\Xi: \cE_S^0|_{V}\xrightarrow[]{\sim} \cE|_{V}.
\]

Working Zariski locally on $S$, we may assume that $S=\on{Spec}(R)$ for a $k$-algebra $R$. Let $x$ be a local parameter around $q$. Let $\mathcal{E}=(E, E_{\bullet}, \nabla)$. We denote $\Delta\coloneqq\nabla_{x\partial_x}$, which is an $R$-linear operator acting on $E$. Since $x\partial_x$ preserves the
ideal $(x)$, $\Delta$ induces an endomorphism $\Delta_q$ of $E_q\coloneqq E|_{S\times q}$, and $\Delta_q=\on{res}_q(\nabla)$. Since $\mathcal{E}$ has vanishing $p$-curvature, we have $\psi_{\nabla}(x\partial_x)=\Delta^p-\Delta=0$. For each $c\in \mathbb{F}_p$, let $P_c(T)$ be the
Lagrange interpolation polynomial
\[
P_c(T)
=\prod_{\substack{a\in\mathbb{F}_p\\ a\neq c}}
\frac{T-a}{c-a}
\in \mathbb{F}_p[T].
\]
Modulo $T^p-T$, these polynomials satisfy
\[
(T-c)P_c(T)=0, \qquad P_c(T)^2=P_c(T), \qquad \sum_{c\in \mathbb{F}_p }P_c(T)=1.
\]
For each $c\in\bF_p$, the idempotent $P_c(\Delta_q)$ satisfies \[\on{Im}P_c(\Delta_q)=\ker(\Delta_q-c)=\ker(\on{res}_q(\nabla)-c).
\]
Hence the eigenspace decomposition of
$\Delta_q$ acting on $E_q$ is given by
\[
E_q=\bigoplus_{c\in\bF_p}\on{Im}P_c(\Delta_q).
\]
Since $\on{res}_q(\nabla)=\Delta_q$ respects the parabolic structure $E_{\bullet}$, the eigenspace decomposition is compatible with the parabolic filtration. It follows that the parabolic weights and filtration type of the induced parabolic structure on each $\on{Im}P_c(\Delta_q)$ are given by
\[
\bigl\{\beta_{\on{op}}^{\widetilde{\evtau}}(i)\mid \tau_i=c\bigr\}
\]
counted with multiplicity. 

After possibly shrinking $S$ to a smaller affine open subset, we choose a basis $\boldsymbol{e}=(e_1, e_2,\dots, e_n)$ of $E_q$ that satisfies the following two properties:
\begin{enumerate}
    \item For each $1\leq i\leq n$, $e_i\in \on{Im}P_{\tau_i}(\Delta_q)$.
    \item For each $c\in\bF_p$ and $1\leq j\leq r$, the images of
    \[
    \bigl\{e_i\mid \tau_i=c,\     \beta_{\on{op}}^{\widetilde{\evtau}}(i)=\beta_j\bigr\}
    \]
    form a basis of the graded quotient $\on{gr}_{\beta_j}\bigl(\on{Im}P_c(\Delta_q)\bigr)$.
\end{enumerate}

Now let $U\subset C$ be an affine open neighborhood of $q$ on which the chosen local parameter $x$ is defined and such that $\frac{dx}{x}$ trivializes $\omega_C(q)|_U$. Since $S\times U$ is affine, for each $1\leq i\leq n$ we may lift
$e_i\in \on{Im}P_{\tau_i}(\Delta_q)$
to an element $\widetilde e_i\in \on{Im}P_{\tau_i}(\Delta|_U).$
We denote the resulting lifts by
$\widetilde{\boldsymbol e}
=(\widetilde e_1,\widetilde e_2,\dots,\widetilde e_n)$. Let $V\subset S\times C$ be the Zariski open locus on which
$(\widetilde e_1,\widetilde e_2,\dots,\widetilde e_n)$ forms a basis of $E$.
Since $(e_1,\dots,e_n)$ is a basis of $E_q$, the open subset $V$ contains
$S\times q$. Since $(\Delta-c)P_c(\Delta)=0$, we have $\Delta(\widetilde{e}_i)=\tau_i \widetilde{e}_i$ for each $1\leq i\leq n$.  
Now we construct the desired isomorphism $\Xi: \cE_S^0|_{V}\xrightarrow[]{\sim} \cE|_{V}$.
Recall that $\mathcal{E}^0=(E^0, E^0_{\bullet},\nabla^0)$ with $E^0=\bigoplus_{i=1}^n \mathcal{O}_C\bigl(-\nu(\tau_i)q\bigr)$. Over the open subset $U\subset C$, for each $1\leq i\leq n$, we set \[
\widetilde e_i^0
\coloneqq
x^{\nu(\tau_i)}
\in
\mathcal O_C\bigl(-\nu(\tau_i)q\bigr)|_U.
\]
Note that $(\widetilde e_1^0, \widetilde e_2^0, \dots, \widetilde e_n^0)$ forms a basis of $E^0|_U$, and $\Delta(\widetilde e_i^0)=\tau_i\widetilde{e}^0_i$ for each $1\leq i\leq n$. Now define
\[
\Xi\colon \cE_S^0|_V \longrightarrow \cE|_V
\]
by setting
\[
\Xi(\widetilde e_i^0)=\widetilde e_i,
\qquad 1\leq i\leq n.
\]
It is immediate from the construction that $\Xi$ is an isomorphism of logarithmic connections \[
(E_S^0,\nabla_S^0)|_V
\xrightarrow{\sim}
(E,\nabla)|_V.
\]
It also respects the parabolic structures, because the parabolic structure on $\cE^0$ (resp. $\cE$) is defined by assigning the parabolic weight $\beta_{\on{op}}^{\widetilde{\evtau}}(i)$ to $\widetilde{e}_i^0|_q$ (resp. $\widetilde{e}_i|_q$) for each $i$.
\end{proof}

Recall that the parabolic weights $\evalpha$ and filtration type $\evm'$ arise from the assignment $\assa$ constructed from $\asstau$ through (\ref{eqn: asstau to assa}). Let $\widetilde{\evalpha}=\bigl(\alpha_{\on{op}}^{\widetilde{\eva}}(1), \alpha_{\on{op}}^{\widetilde{\eva}}(2), \dots, \alpha_{\on{op}}^{\widetilde{\eva}}(n)\bigr)$. Then we have
\begin{lemma}\label{lemma: connection auto and parabolic group scheme}
There is a canonical isomorphism of sheaves of groups on $C'$
\begin{equation}\label{eqn: identifying auto with parabolic}
 (\Fr_C)_*{\Aut}(\cE^0)
\xrightarrow{\sim}
\cG_{\widetilde{\evalpha}}.   
\end{equation}
\end{lemma}

\begin{proof}
 Away from $q'$, the sheaf of groups $(\Fr_C)_*{\Aut}(\cE^0)$ is canonically isomorphic to $\on{GL}_n$ by Cartier descent, see \cite[Theorem~5.1]{Katz1970}. Therefore we have an isomorphism
 \[
(\Fr_C)_*{\Aut}(\cE^0)|_{C'\setminus q'}
\xrightarrow{\sim}
\cG_{\widetilde{\evalpha}}|_{C'\setminus q'}.
 \]
 The fact that this isomorphism extends over $q'$ essentially follows from the computation in  \cite[Lemma~4.14]{li-sun2021tame}, which we spell out in detail for the reader's convenience. Let $D_q=\on{Spec}k[[x]]$ be the formal disc around $q\in C$. An element
\[
g\in \Aut(\cE^0)|_{D_q}
\subset \on{GL}_n\bigl(k((x))\bigr)
\]
is characterized by the following two conditions:
\begin{enumerate}
    \item \label{cond:connection}$g$ preserves the connection, i.e.
    \[dg=gd,
    \]
    \item \label{cond:lattice}$g$ preserves the lattice underlying $E^0|_{D_q}$ and the parabolic structure at $q$, i.e.
    \[x^{-\nu(\widetilde{\evtau})}gx^{\nu(\widetilde{\evtau})}
    \in \cG_{\widetilde{\evbeta}}(k[[x]]), \qquad\text{where } \widetilde{\evbeta}=\bigl(\beta_{\on{op}}^{\widetilde{\evtau}}(1), \beta_{\on{op}}^{\widetilde{\evtau}}(2), \dots, \beta_{\on{op}}^{\widetilde{\evtau}}(n)\bigr).
    \]
\end{enumerate}
We write $g=(g_{ij})$, then condition (\ref{cond:connection}) is equivalent to 
\[
g_{ij}=\sum_{\ell} a_{ij}^{(\ell)}x^{\ell p}.
\]
Note that $x^{-\nu(\widetilde{\evtau})}gx^{\nu(\widetilde{\evtau})}=(g'_{ij})$, where $g'_{ij}=\sum_{\ell} a_{ij}^{(\ell)}x^{\ell p-\nu({\tau_i})+\nu({\tau_j})}$. 
Therefore conditions (\ref{cond:connection}) and (\ref{cond:lattice}) combine to
\begin{equation}\label{eqn: connection eq}
\ell p-\nu(\tau_i)+\nu(\tau_j)
\geq
\left\lceil
\beta_{\on{op}}^{\widetilde{\evtau}}(j)
-
\beta_{\on{op}}^{\widetilde{\evtau}}(i)
\right\rceil \quad\text{whenever }a_{ij}^{(\ell)}\neq 0, \qquad \det(g)\in k[[x^p]]^\times.
\end{equation}
It follows that (\ref{eqn: connection eq}) is equivalent to $g
    \in \cG_{\widetilde{\evalpha}}(k[[x^p]])$.
The same calculation works over an arbitrary
$k$-algebra $R$.
\end{proof}

Before proving \Cref{thm: identifying with parabolic higgs},
we first establish an isomorphism between the stacks
obtained from the two sides of \Cref{thm: identifying with parabolic higgs} by forgetting the datum $\psi$ and
the Higgs field, respectively.

\begin{prop}\label{prop: logarithmic Cartier descent}
For any $\evtau\in \mathfrak{t}(\bF_p)/W$, parabolic weights $\evbeta$, and assignment $\asstau$, there exists an isomorphism
\begin{equation}  \label{eqn: equivalence between connection and parabolic bundle}\con_0^{\asstau}(C)\cong \Bun^{\evalpha}(C', \evm'), 
\end{equation}   
where the parabolic weights $\evalpha$ and filtration type $\evm'$ arise from the assignment $\assa$ constructed from $\asstau$ through (\ref{eqn: asstau to assa}). 
\end{prop}
\begin{proof}
For a sheaf of groups $\cG$ on $C$, we denote by $\Bun_{\cG}(C)$ the moduli stack of $\cG$-torsors on $C$. By Lemma~\ref{lemma: connection auto and parabolic group scheme} and Remark~\ref{remark: parabolic bundle and parabolic group torsor}, we have 
\[
\Bun^{\evalpha}(C', \evm')\cong \Bun_{(\Fr_C)_*{\Aut}(\cE^0)}(C').
\]
Note that $\Fr_C\colon C\to C'$ induces an equivalence of \'etale sites, therefore we have
\[
\Bun_{(\Fr_C)_*{\Aut}(\cE^0)}(C')\cong \Bun_{\Aut(\cE^0)}(C). 
\]
It follows that the desired (\ref{eqn: equivalence between connection and parabolic bundle}) is equivalent to
\[
\con_0^{\asstau}(C)\cong \Bun_{\Aut(\cE^0)}(C). 
\]
Applying Lemma \ref{lemma: local isomorphism}, we define 
\[
\begin{aligned}
\mathsf{F}\colon
\con_0^{\asstau}(C)
&\longrightarrow
\Bun_{\Aut(\cE^0)}(C),\\
\cE
&\longmapsto
\Isom(\cE^0,\cE).
\end{aligned}
\]
In the reverse direction, we define
\[
\begin{aligned}
\mathsf{G}\colon
\Bun_{\Aut(\cE^0)}(C)
&\longrightarrow
\con_0^{\asstau}(C),\\
\cP
&\longmapsto
\cP\times^{\Aut(\cE^0)}\cE^0,
\end{aligned}
\]
with connection and parabolic structure inherited from $\cE^0$. Consider the natural morphism 
\[
\varepsilon_{\cE}: \mathsf{G}\circ \mathsf{F}(\cE)=\Isom(\cE^0,\cE)\times^{\Aut(\cE^0)}\cE^0\longrightarrow \cE
\]
induced by evaluation, sending $(\varphi, e)$ to $\varphi(e)$. It follows from Lemma \ref{lemma: local isomorphism} that $\varepsilon_{\cE}$ is an isomorphism of parabolic connections. Moreover, under these identifications, the maps on Hom-sets induced by
$\mathsf{F}$ and $\mathsf{G}$ are inverse to each other. Hence $\mathsf{F}$ and $\mathsf{G}$ are isomorphisms of stacks that are inverse to each other.    
\end{proof}
\begin{proof}[Proof of \Cref{thm: identifying with parabolic higgs}]

We now obtain the desired isomorphism 
\[\mathfrak{X}^{\asstau}\cong \hig^{\evalpha}(C', \evm') 
\]
from the isomorphism (\ref{eqn: equivalence between connection and parabolic bundle}) in \Cref{prop: logarithmic Cartier descent}. Recall that $\con_0^{\asstau}(C)$ is obtained from $\mathfrak{X}^{\asstau}$ by forgetting the datum $\psi$, which is an $\mathcal{O}_C$-linear morphism $\psi: E\to E\otimes \Fr_C^{*}\bigl(\omega_{C'}(q')\bigr)$ satisfying the two conditions (\ref{cond: horizontal}) and (\ref{cond: psi respects flag}) described earlier. For any $\cE=(E, E_{\bullet}, \nabla)\in \con_0^{\asstau}(C)$ and any line bundle $\cL$ on $C'$, we have 
\[
\cE\otimes \Fr_C^*\cL\coloneqq\left(
E\otimes\Fr_C^*\cL,\,
E_{\bullet}\otimes\Fr_C^*\cL
,\,\nabla\otimes 1+1\otimes\nabla_{\mathrm{can}}\right)\in \con_0^{\asstau}(C)
\]
where $\nabla_{\mathrm{can}}$ denotes the canonical connection on $\Fr_C^*\cL$. Note that the isomorphism (\ref{eqn: identifying auto with parabolic})
\[
(\Fr_C)_*\Aut(\cE^0)\cong\cG_{\widetilde{\evalpha}}
\]
in \Cref{lemma: connection auto and parabolic group scheme}
identifies the central scalar automorphisms on both sides. Therefore under the isomorphism~(\ref{eqn: equivalence between connection and parabolic bundle})
\[\con_0^{\asstau}(C)\cong \Bun^{\evalpha}(C', \evm'),\]
if $\cE$ corresponds to $(F, F_{\bullet})$, then $\cE\otimes\Fr_C^*\cL$ corresponds to $(F\otimes \cL, F_{\bullet}\otimes \cL)$. Now set $\cL=\omega_{C'}(q')$. Then the two conditions on $\psi$ are equivalent to requiring that
\[
\psi\colon \cE\longrightarrow
\cE\otimes \Fr_C^*\bigl(\omega_{C'}(q')\bigr)
\]
be a morphism of parabolic connections. Under the isomorphism~(\ref{eqn: equivalence between connection and parabolic bundle}), $\psi$ corresponds to a morphism $F\to F\otimes \omega_{C'}(q')$ that respects the parabolic structures $F_{\bullet}$ and $F_{\bullet}\otimes \omega_{C'}(q')$, which is precisely a parabolic Higgs field. The $\cP$-equivariance follows from the isomorphism \eqref{eqn: identifying auto with parabolic}
\[(\Fr_C)_*\Aut(\cE^0)\cong\cG_{\widetilde{\evalpha}}
\]
and the induced identification of their Lie algebras, which matches $\psi$ with the parabolic Higgs field on the other side after locally trivializing $\omega_{C'}(q')$. Consequently,
tensoring the spectral sheaf of $\psi$ by $W^*L$
corresponds to tensoring the spectral sheaf of the
associated parabolic Higgs field by $L$.
\end{proof}

\subsubsection{Proof of \Cref{thm: stacky main para}}\;

Now we fix $\eva\in \mathfrak{t}(k)/W$. Let $\mathfrak{X}_{\eva}^{\asstau}\coloneqq \mathfrak{X}^{\asstau}\times_{A}A_{\eva}$ be the restriction of $\mathfrak{X}^{\asstau}$ to $A_{
\eva}$. The stack $\mathfrak{X}_{\eva}^{\asstau}$ admits a further decomposition into a disjoint union, which we now describe. Let $(E, E_{\bullet},\nabla, \psi)\in \mathfrak{X}_{\eva}^{\asstau}$. Since $\mathfrak{X}_{\eva}^{\asstau}$ lies over $A_{
\eva}$, the residue $\on{res}_q(\psi)$ has eigenvalues $\eva$. Since $\psi$ is horizontal with respect to the canonical logarithmic connection constructed from $\nabla$, $\on{res}_q(\psi)$ and $\on{res}_q(\nabla)$ are commuting matrices. By definition, both $\on{res}_q(\psi)$ and $\on{res}_q(\nabla)$ respect the flag structure $E_{\bullet}$. It follows that the action of $\on{res}_q(\psi)$ and $\on{res}_q(\nabla)$ induces a decomposition into generalized joint eigenspaces
\[E_q=\bigoplus_{(a,\tau)\in \eva\times \evtau}V^{a, \tau},
\]
and the parabolic structure $(E_q, E_{\bullet}, \evbeta)$ induces a parabolic structure on each $V^{a, \tau}$. Let $\evb^{\mathrm{reg}}\in \mathfrak{t}(k)/W$ be the unique regular tuple satisfying $\on{AS}(\evb^{\mathrm{reg}})=\eva$ and $\nu(\evb^{\mathrm{reg}})=0$. We define $\evb\in \mathfrak{t}(k)/W$ to be the unordered $n$-tuple in which
$b+\tau$ occurs with multiplicity $\dim_k V^{\on{AS}(b),\tau}$
for each $(b,\tau)\in \evb^{\mathrm{reg}}\times \evtau$. The parabolic structures on all $V^{\on{AS}(b),\tau}$ determine an
assignment $\assb$ as follows. For each parabolic weight $\beta_j\in \evbeta$ and each $(b,\tau)\in \evb^{\mathrm{reg}}\times \evtau$, we assign
the weight $\beta_j$ to $\dim_k \on{gr}_{\beta_j}
V^{\on{AS}(b),\tau}$ copies of $b+\tau$. In this way, we obtain an assignment
$\assb$ of the same filtration type $\evm$ as $\asstau$. 
\begin{defn}\label{defn: assb substack of p-curvature zero}
  We define $\mathfrak{X}_{\eva}^{\assb}\subseteq \mathfrak{X}_{\eva}^{\asstau}$ to be the substack characterized by the actions of $\on{res}_q(\psi)$ and $\on{res}_q(\nabla)$ on $(E_q, E_{\bullet})$ giving rise to the assignment $\assb$ according to the procedure described above.  
\end{defn}
Then we have
\[
\mathfrak{X}_{\eva}^{\asstau}=\coprod\mathfrak{X}_{\eva}^{\assb},
\]
where the disjoint union ranges over assignments $\assb$ of filtration type $\evm$ whose $\evb$-component satisfies $\on{AS}(\evb)=\eva$ and $\evb-_{\eva}\evb^{\mathrm{reg}}=\evtau$. Here $(\evb-_{\eva}\evb^{\mathrm{reg}})\in \mathfrak{t}(k)/W$ is obtained from $\evb$ by replacing each occurrence of $b$ in $\evb$ with $b-b'$, where $b'$ is the unique element of $\evb^{\mathrm{reg}}$ satisfying $\on{AS}(b)=\on{AS}(b')$. 

Restricting the isomorphism (\ref{eqn: Li-Sun naH}) 
\[
\cH\times_A^{\mathcal P^{\text{ext}}}\mathfrak{X}^{\evbeta, \evm}\xrightarrow{\sim} \con^{\evbeta}(C, \evm)
\]
to smaller substacks, we obtain the following
\begin{prop}\label{prop: Li-Sun fixed assignment}
    There exists an isomorphism 
\begin{equation}  \cH_{\eva}^{\good}\times_{A_{\eva}}^{\mathcal P_{\eva}^o}\mathfrak{X}_{\eva}^{\assb}\xrightarrow{\sim} \con^{\assb}(C).
\end{equation}
\end{prop}
\begin{proof}
    Let $\cP'\subset \cP^{\on{ext}}$ be the substack characterized by $\on{res}_q(\nabla)=0$. By definition, $\cP'$ identifies with the substack of $\mathfrak{X}^{\{0\},(n),\boldsymbol{0}}$ characterized by the spectral sheaf of $\psi$ being invertible. Here $\boldsymbol{0}$ denotes the zero tuple of residue
eigenvalues, $\{0\}$ denotes the trivial parabolic weight, and $(n)$ denotes the trivial filtration type. By \Cref{thm: identifying with parabolic higgs} and the BNR correspondence, we have $\cP'\cong\cP(S/A)=\cP$, and its actions on $\cH$ and $\mathfrak{X}^{\evbeta, \evm}$ agree with the ones considered in Section \ref{sec: torsor conn} and \Cref{thm: identifying with parabolic higgs}, respectively. Let $\evb^{\mathrm{reg}}\in \mathfrak{t}(k)/W$ be the unique regular tuple satisfying $\on{AS}(\evb^{\mathrm{reg}})=\eva$ and $\nu(\evb^{\mathrm{reg}})=0$. Let $\mathfrak{X}_{\eva}^{\evbeta,\evm}\coloneqq \mathfrak{X}^{\evbeta,\evm}\times_{A}A_{\eva}$. Then the isomorphism (\ref{eqn: Li-Sun naH}) restricts to an isomorphism
  \begin{equation}\label{eqn: Li-Sun very good over a}   \cH_{\evb^{\on{reg}}}^{\good}\times_{A_{\eva}}^{\mathcal P_{\eva}^o}\mathfrak{X}_{\eva}^{\evbeta, \evm}\xrightarrow{\sim} 
\coprod_{\substack{\evb\in\mathfrak t(k)/W\\
\on{AS}(\evb)=\eva}}
\con^{\evbeta}_{\evb}(C,\evm),
  \end{equation}   
  where $\con^{\evbeta}_{\evb}(C,\evm)\subset \con^{\evbeta}(C,\evm)$ is the substack characterized by the multiset of eigenvalues of $\on{res}_q(\nabla)$ being equal to $\evb$. Recall that $\cH_{\eva}^{\good}$ is identified with $\cH_{\evb^{\on{reg}}}^{\good}$ through \Cref{prop: heck transform for good connections}. Note that 
    \[   \mathfrak{X}_{\eva}^{\evbeta, \evm}=\coprod\mathfrak{X}_{\eva}^{\assb},
    \]
    where the disjoint union ranges over all $\evb\in \mathfrak{t}(k)/W$ satisfying $\on{AS}(\evb)=\eva$ and all assignments $\assb$ of filtration type $\evm$. We also have 
\[
\con^{\evbeta}_{\evb}(C,\evm)=\coprod\con^{\assb}(C)
\]
where the disjoint union ranges over all assignments $\assb$ of filtration type $\evm$. Now we check that the isomorphism (\ref{eqn: Li-Sun very good over a}) matches $\mathfrak{X}_{\eva}^{\assb}$ with $\con^{\assb}(C)$. Let $\cE=(E, \nabla_E)\in \cH^{\good}_{\evb^{\on{reg}}}$ such that $h_{\on{dR}}(\cE)=a\in A_{\eva}$, and let $L_{\cE}$ be the spectral sheaf of its $p$-curvature $\psi_{\nabla_E}$. Let $\cM=(M, M_{\bullet},\nabla_M, \psi)\in \mathfrak{X}^{\assb}_{\eva}$ that maps to the same $a\in A_{\eva}$, and let $\cF_{\psi}$ be the spectral sheaf of $\psi$. Recall that the isomorphism (\ref{eqn: Li-Sun very good over a}) maps $(\cE, \cM)$ to 
\[
\left(\pi_*^{\Fr}(L_\cE\otimes_{S_a^{\Fr}} \cF_{\psi}), \bigl(\pi_*^{\Fr}(L_\cE\otimes_{S_a^{\Fr}} \cF_{\psi})\bigr)_{\bullet}, \nabla_{\otimes}\right),
\]
where $\nabla_{\otimes}$ is the tensor product connection, and the parabolic structure $\bigl(\pi_*^{\Fr}(L_\cE\otimes_{S_a^{\Fr}} \cF_{\psi})\bigr)_{\bullet}$ is induced from the parabolic structure on $\pi_*^{\Fr}\cF_{\psi}$. Recall that the spectral curve $S_a^{\Fr}$ decomposes as 
\[S_a^{\Fr}\times_C D_q=\coprod_{j=1}^s S_{a,a_j}^{\Fr} 
\]
near $q$, and the spectral sheaf $\cF_{\psi}$ decomposes accordingly, with summands corresponding to different eigenvalues of $\on{res}_q(\psi)$. Let $b_j\in \evb^{\on{reg}}$ be the unique element such that $\on{AS}(b_j)=a_j$. By definition, $\on{res}_q(\nabla_E)$ has eigenvalues equal to $b_j$ when restricted to $\pi^{\Fr}_*(L_{\cE}|_{S_{a,a_j}^{\Fr}})$. Tensoring with $L_{\cE}$ over $S_{a,a_j}^{\Fr}$ adds $b_j$ to the eigenvalues of $\on{res}_q(\nabla_M)$, while transporting the parabolic structure on $\pi_*^{\Fr}\cF_{\psi}$ accordingly. This modification of eigenvalues and parabolic structure agrees exactly with the modification from $\asstau$ to $\assb$ used in \Cref{defn: assb substack of p-curvature zero} of $\mathfrak{X}_{\eva}^{\assb}$.
\end{proof}
Note that when restricted to $A_{\eva}\subset A$, the right-hand side of the isomorphism in \Cref{thm: identifying with parabolic higgs} decomposes as 
\[
\hig^{\evalpha}(C', \evm')\times_{A}A_{\eva}=\coprod \hig^{\assa}(C'),
\]
where the disjoint union ranges over all assignments $\assa$ of filtration type $\evm'$. We have

\begin{prop}\label{prop: refinement of cartier descent}
    Let $\assa$ and $\assb$ be two assignments related through (\ref{eqn: assa from assb}). The isomorphism in \Cref{thm: identifying with parabolic higgs} restricts to an isomorphism
\begin{equation} 
\mathfrak{X}_{\eva}^{\assb}\cong \hig^{\assa}(C').
\end{equation}
\end{prop}
\begin{proof}
  Let $\cE^0=(E^0, E^0_{\bullet},\nabla^0)\in \con_0^{\asstau}(C)(k)$ be the model object defined in \Cref{sec: proof of identifying with parabolic higgs}. Let $D_q=\on{Spec}k[[x]]$ be the formal disc around $q\in C$. The same computation as  in the proof of \Cref{lemma: connection auto and parabolic group scheme} implies that
  \[\phi\in \End(\cE^0)(D_q)\subset \on{Lie}(\on{GL}_n)\bigl(k((x))\bigr) \Longleftrightarrow \phi\in \on{Lie}(\cG_{\widetilde{\evalpha}})\bigl(k[[x^p]]\bigr)\Longrightarrow \psi\coloneqq x^{-\nu(\widetilde{\evtau})}\phi x^{\nu(\widetilde{\evtau})}
    \in \on{Lie}(\cG_{\widetilde{\evbeta}})\bigl(k[[x]]\bigr),\]
  where $\widetilde{\evalpha}=\bigl(\alpha_{\on{op}}^{\widetilde{\eva}}(1), \alpha_{\on{op}}^{\widetilde{\eva}}(2), \dots, \alpha_{\on{op}}^{\widetilde{\eva}}(n)\bigr)$ as defined through (\ref{eqn: asstau to assa}), and $\widetilde{\evbeta}=\bigl(\beta_{\on{op}}^{\widetilde{\evtau}}(1), \beta_{\on{op}}^{\widetilde{\evtau}}(2), \dots, \beta_{\on{op}}^{\widetilde{\evtau}}(n)\bigr)$.
Let 
\[
l_{\widetilde{\evalpha}}: \on{Lie}(\cG_{\widetilde{\evalpha}})\bigl(k[[x^p]]\bigr)\to \on{Lie}(L_{\widetilde{\evalpha}})(k)
\]
be the morphism on Lie algebras induced by the quotient map from a parahoric group scheme to the corresponding Levi factor. Let $\phi\in \End(\cE^0)(D_q)$ such that the multiset of eigenvalues of $\phi_q$ acting on $E^0_q$ is given by $\eva$. To obtain the desired isomorphism, it is enough to show that if the actions of $\phi_q$ and $\on{res}_q(\nabla^0)$ on $(E_q^0, E^0_{\bullet})$ give rise to $\assb$ through the procedure used in Definition \ref{defn: assb substack of p-curvature zero}, then the block diagonal matrix $l_{\widetilde{\evalpha}}(\phi)$ has eigenvalues prescribed by $\assa$. Note that $E^0$ is trivialized through the local section $x^{\nu(\widetilde{\evtau})}$, and the action of $\phi_q$ becomes $\psi(0)$. A generalized joint eigenvector of $(\psi(0),\on{res}_q(\nabla^0))$ with eigenvalues $(a,\tau)\in \eva\times\evtau$ and parabolic weight $\beta \in \evbeta$ corresponds to a generalized eigenvector of $\phi(0)$ with eigenvalue $a$ and parabolic weight $\alpha=(\beta-\nu(\tau))/p$. If $b=b'+\tau$ with $b'\in\evb^{\mathrm{reg}}$
and $\tau\in\bF_p$, then $\nu(b)=\nu(\tau)$. The desired statement follows from simple unwinding of notations. 
\end{proof}
\Cref{thm: stacky main para} follows immediately by combining \Cref{prop: Li-Sun fixed assignment} and \Cref{prop: refinement of cartier descent}. 

\subsection{Proof of the main theorem}
We write ${\cC}onn_0^{\asstau,\mathrm{ss}}(C)\coloneqq {\cC}onn_0^{\asstau}(C)\cap {\cC}onn^{\evbeta\text{-}\mathrm{ss}}(C,\evm)$. We start by showing 
\begin{prop}\label{prop: semistable Cartier descent}
 The isomorphism 
\begin{equation}\label{eqn: proving cartier descent match semistable}
\con_0^{\asstau}(C)\cong \Bun^{\evalpha}(C', \evm')
\end{equation}
in Proposition~\ref{prop: logarithmic Cartier descent} restricts to an isomorphism
\[
\con_0^{\asstau, \mathrm{ss}}(C)\cong \Bun^{\evalpha\text{-}\mathrm{ss}}(C', \evm').
\]
\end{prop}
\begin{proof}
    Let $\cE=(E, E_{\bullet},\nabla)\in \con_0^{\asstau}(C)$ and let $\cF=(F, F_{\bullet})\in \Bun^{\evalpha}(C', \evm')$ be its image under the isomorphism (\ref{eqn: proving cartier descent match semistable}). Note that for any vector bundle $E$ over $C$, any vector subbundle $E'\subseteq E$ is uniquely determined by its restriction $E'|_{U}$ to $U=C\setminus q$. Indeed, writing $j:U\hookrightarrow C$ for the inclusion, we have
\[
E'=E\cap j_*(E'|_U)
\]
inside $j_*(E|_U)$. Over $U$ and $U'$, the correspondence (\ref{eqn: proving cartier descent match semistable}) is given by Cartier descent, i.e.
\[
F|_{U'}\cong (\Fr_C)_*\ker(\nabla|_U:E|_U\to E|_U\otimes\omega_U), \quad (E|_U, \nabla|_{U})\cong (\Fr_C^*(F|_{U'}), \nabla_{\on{can}}).
\]
Cartier descent also induces a correspondence between $\nabla|_U$-invariant subbundles of $E|_U$ and subbundles of $F|_{U'}$. Let $E_U'\subseteq E|_U$ be the $\nabla|_U$-invariant subbundle corresponding to $F'_{U'}\subseteq F|_{U'}$. Let $E'=E\cap j_*(E'_U)$ and $F'=F\cap j'_*(F'_{U'})$. These are the unique saturated extensions of $E'_U$ and $F'_{U'}$, respectively. Let $(E', E_{\bullet}', \evbeta')$ (resp. $(F', F_{\bullet}', \evalpha')$) be the parabolic subbundle of $(E,E_\bullet,\evbeta)$ (resp. $(F,F_\bullet,\evalpha)$)
with the induced parabolic structure. For the desired matching of semistable substacks through (\ref{eqn: proving cartier descent match semistable}), it is enough to show that for any choice of $E'_U$ and $F'_{U'}$, we have 
\begin{equation}\label{eqn: matching parabolic degree}
\on{pardeg}(E', E_{\bullet}', \evbeta')=p\on{pardeg}(F', F_{\bullet}', \evalpha'). 
\end{equation}

Recall that the degree of a vector bundle $E$ on $C$ is equal to the degree of its determinant line bundle $L=\on{det}(E)$, which can be computed as the sum of the orders of vanishing of any nonzero rational section of $L$ at all points of $C$, with poles counted negatively. Since $E'_U\cong \Fr^*_CF'_{U'}$, the contributions over $U$ and $U'$ to the respective parabolic degrees satisfy the multiplication-by-$p$ formula (\ref{eqn: matching parabolic degree}). 

Now we compute the contributions at $q$ and $q'$, which can be done by working with the model object $\cE^0=(E^0, E^0_{\bullet},\nabla^0, \evbeta)$ defined in Section \ref{sec: proof of identifying with parabolic higgs}, as well as the model object $\cF^0=(F^0, F^0_{\bullet}, \evalpha)$ associated to $\widetilde{\evalpha}$ as in \Cref{remark: parabolic bundle and parabolic group torsor}. Let $D=\on{Spec}(k[[x]])$ be the formal disc around $q\in C$. Let $e_1, e_2, \dots, e_n$ be the standard basis of $F^0|_{D'}\cong\bigoplus_{i=1}^n \mathcal{O}_{D'}$. We equip $\{1,2,\dots,n\}$ with a total order by setting
\[
i\prec j \quad\text{if}\quad \alpha_{\on{op}}^{\widetilde{\eva}}(i)<\alpha_{\on{op}}^{\widetilde{\eva}}(j),
\]
and when $\alpha_{\on{op}}^{\widetilde{\eva}}(i)=\alpha_{\on{op}}^{\widetilde{\eva}}(j)$, we break ties using the usual ordering of the indices. By choosing an arbitrary basis of $F'_{q'}\subseteq k^n$ and applying Gaussian elimination, we obtain a basis $v_1, v_2, \dots, v_{n'}$ of the form
\[
v_j=e_{i_j}+\sum_{i_j\prec i}c_{ji}e_i,
\qquad i_1\prec i_2\prec\cdots\prec i_{n'}.
\]
Consequently, the parabolic weights (counted with multiplicity) of $F'$ at $q'$ are
$\alpha_{\on{op}}^{\widetilde{\eva}}(i_1),\alpha_{\on{op}}^{\widetilde{\eva}}(i_2), \dots,\alpha_{\on{op}}^{\widetilde{\eva}}(i_{n'})$, and their contribution
to the parabolic degree is $\sum_{j=1}^{n'}\alpha_{\on{op}}^{\widetilde{\eva}}(i_j)$. Each $v_j$ lifts to a section $\widetilde{v}_j$ of $F'|_{D'}$, which is of the form
\[
\widetilde{v}_j=e_{i_j}+\sum_{i_j\prec i}c_{ji}e_i
+x^p\sum_i a_{ji}(x^p)e_i.
\]
$\widetilde{v}_1, \widetilde{v}_2, \dots, \widetilde{v}_{n'}$ form a basis of $F'|_{D'}$.

By abuse of notation, we also denote by $e_1, e_2, \dots, e_n$ the standard basis of \[\Fr^*_C(F^0|_{D'})=\displaystyle\bigoplus_{i=1}^n \mathcal{O}_D\subseteq E^0|_D=\displaystyle\bigoplus_{i=1}^n \mathcal{O}_D\bigl(-\nu(\tau_i)q\bigr).
\]
Set $f_i=x^{\nu{(\tau_i)}}e_i$, then $E^0|_D$ is trivialized as $E^0|_D=\bigoplus_{i=1}^n \mathcal{O}_Df_i$. Now set 
$\widetilde{u}_j\coloneqq x^{\nu(\tau_{i_j})}\widetilde{v}_j$. Then
\[
\widetilde{u}_j=f_{i_j}+\sum_{i_j\prec i}c_{ji}x^{\nu(\tau_{i_j})-\nu(\tau_i)}f_i
+\sum_i x^{p+\nu(\tau_{i_j})-\nu(\tau_i)}a_{ji}(x^{p})f_i.
\]
Note that
\begin{equation}\label{eqn: index ordering through parabolic weights}
i\prec j \implies \alpha_{\on{op}}^{\widetilde{\eva}}(i)\leq\alpha_{\on{op}}^{\widetilde{\eva}}(j)
\implies
\begin{cases}
\nu(\tau_i)>\nu(\tau_j), \quad\text{or}\\
\nu(\tau_i)=\nu(\tau_j)\text{ and }\beta_{\on{op}}^{\widetilde{\evtau}}(i)\leq\beta_{\on{op}}^{\widetilde{\evtau}}(j).  
\end{cases}
\end{equation}
Also note that since the image of $\nu$ lies in $\{-p+1, \dots, -1,0\}$, we have $p+\nu(\tau_{j})-\nu(\tau_i)>0$ for any indices $i,j$. Therefore $\widetilde{u}_j$ is a section of $E^0|_D$ for each $1\leq j\leq n'$. Furthermore, $\widetilde{u}_1(0), \widetilde{u}_2(0),\dots, \widetilde{u}_{n'}(0)$ are linearly independent, and $\widetilde{u}_1|_{D^{\times}}, \widetilde{u}_2|_{D^{\times}},\dots, \widetilde{u}_{n'}|_{D^{\times}}$ form a basis of $E'|_{D^{\times}}$. By uniqueness of saturated extension, $\widetilde{u}_1, \widetilde{u}_2,\dots, \widetilde{u}_{n'}$ form a basis of $E'|_D$. It follows that the contribution to $\on{deg}(E')$ at $q$ is $\sum_{j=1}^{n'}(-\nu(\tau_{i_j}))$. Note that \[
\widetilde u_j(0)
=f_{i_j}
+\sum_{\substack{i_j\prec i\\
\nu(\tau_i)=\nu(\tau_{i_j})}}
c_{ji}f_i \qquad \text{for each }1\leq j\leq n'.
\]
By (\ref{eqn: index ordering through parabolic weights}), the parabolic weights (counted with multiplicity) of $E'$ at $q$ are $\beta_{\on{op}}^{\widetilde{\evtau}}(i_1), \beta_{\on{op}}^{\widetilde{\evtau}}(i_2),\dots, \beta_{\on{op}}^{\widetilde{\evtau}}(i_{n'})$. The desired equality
\[
\sum_{j=1}^{n'}\beta_{\on{op}}^{\widetilde{\evtau}}(i_{j})+\sum_{j=1}^{n'}(-\nu(\tau_{i_j}))=p\sum_{j=1}^{n'}\alpha_{\on{op}}^{\widetilde{\eva}}(i_j)
\]
follows directly from the definition (\ref{eqn: asstau to assa}).

\end{proof}

\begin{proof}[Proof of \Cref{thm: main para}]
Let $\mathfrak{X}_{\eva}^{\assb,\mathrm{ss}}\subset\mathfrak{X}_{\eva}^{\assb}$ be the substack characterized by the semistability inequality for any subbundle preserved by both the connection $\nabla$ and the twisted Higgs field $\psi$. We start by showing that the isomorphism 
\begin{equation}\label{eqn: stachy correspondence in proof of main theorem}
\cH_{\eva}^{\good}\times_{A_{\eva}}^{\mathcal P_{\eva}^o}\mathfrak{X}_{\eva}^{\assb}\xrightarrow{\sim} \con^{\assb}(C)
\end{equation}
in \Cref{prop: Li-Sun fixed assignment} restricts to an isomorphism
\begin{equation}\label{eqn: restricted Li-Sun semistable}
\cH_{\eva}^{\good}\times_{A_{\eva}}^{\mathcal P_{\eva}^o}\mathfrak{X}_{\eva}^{\assb, \mathrm{ss}}\xrightarrow{\sim} \con^{\assb,\mathrm{ss}}(C).
\end{equation}

Let $\evb^{\mathrm{reg}}\in \mathfrak{t}(k)/W$ be the unique regular tuple satisfying $\on{AS}(\evb^{\mathrm{reg}})=\eva$ and $\nu(\evb^{\mathrm{reg}})=0$. We identify $\cH_{\eva}^{\good}$ with $\cH_{\evb^{\on{reg}}}^{\good}$ through \Cref{prop: heck transform for good connections}. Let $\cE=(E,\nabla_E)\in \cH_{\evb^{\on{reg}}}^{\good}$ such that $h_{\on{dR}}(\cE)=a\in A_{\eva}$, and let $L_{\cE}$ be the spectral sheaf of its $p$-curvature $\psi_{\nabla_E}$. Let $\cM=(M, M_{\bullet},\nabla_M, \psi)\in \mathfrak{X}^{\assb}_{\eva}$ that maps to the same $a\in A_{\eva}$, and let $\cF_{\psi}$ be the spectral sheaf of $\psi$. Recall that the isomorphism (\ref{eqn: stachy correspondence in proof of main theorem}) maps $(\cE, \cM)$ to 
\[
\cN=(N, N_{\bullet}, \nabla_N)\coloneqq \left(\pi_*^{\Fr}(L_\cE\otimes_{S_a^{\Fr}} \cF_{\psi}), \bigl(\pi_*^{\Fr}(L_\cE\otimes_{S_a^{\Fr}} \cF_{\psi})\bigr)_{\bullet}, \nabla_{\otimes}\right), 
\]
see (\ref{eqn: stacky parabolic naH without torsor}). Let $M'\subseteq M$ be a subbundle preserved by $\nabla_M$ and $\psi$. Let $Y_a^{\Fr}\subseteq S_a^{\Fr}$ be the spectral subcurve defined as the vanishing scheme of the characteristic polynomial of $\psi|_{M'}$, and let $\cF'_{\psi}$ be the spectral sheaf on $Y_a^{\Fr}$ corresponding to $M'$. Set $N'\coloneqq \pi_*^{\Fr}(L_\cE|_{Y_a^{\Fr}}\otimes_{Y_a^{\Fr}} \cF'_{\psi})$. Then $N'\subseteq N$ is a subbundle preserved by $\nabla_N$, and every subbundle of $N$ preserved by $\nabla_N$ arises this way from some subbundle $M'\subseteq M$ preserved by $\nabla_M$ and $\psi$. Note that if $S^{\Fr}_{a}$ decomposes into integral components as
\[S^{\Fr}_{a}=\sum_{t=1}^N n_t S^{\Fr}_{a,t},
\]
then 
\[
Y^{\Fr}_{a}=\sum_{t=1}^N n'_t S^{\Fr}_{a,t}, \qquad 0\leq n'_t\leq n_t.
\]
Since $\nu(\evb^{\mathrm{reg}})=0$,
\Cref{defn: good component} implies that $L_{\cE}$ has degree zero on each component $S^{\Fr}_{a,t}$. Hence, by
\cite[Proof of Theorem 4.7, (22)]{dCGZ}, we have 
\[
\chi(\cF'_{\psi})=\chi(L_\cE|_{Y_a^{\Fr}}\otimes_{Y_a^{\Fr}} \cF'_{\psi}), 
\]
thus
\[
\chi(\pi_*^{\Fr}\cF'_{\psi})=\chi(\pi_*^{\Fr}(L_\cE|_{Y_a^{\Fr}}\otimes_{Y_a^{\Fr}} \cF'_{\psi})),
\]
which implies 
\[
\on{deg}(M')=\on{deg}(N').
\]
Recall that the flag structure $N_{\bullet}=\bigl(\pi_*^{\Fr}(L_\cE\otimes_{S_a^{\Fr}} \cF_{\psi})\bigr)_{\bullet}$ is induced from the flag structure  $M_{\bullet}=(\pi_*^{\Fr}\cF_{\psi})_{\bullet}$. It follows that the induced flag structure on $N'_q$ arises in the same fashion from the induced flag structure on $M'_q$. Thus they have the same multiset of parabolic weights counted with multiplicity. It follows that
\begin{equation}\label{eqn: matching degree of connections}
\on{pardeg}(M', M'_{\bullet})=\on{pardeg}(N', N'_{\bullet}),
\end{equation}
which implies the desired isomorphism 
(\ref{eqn: restricted Li-Sun semistable}). 

By \Cref{prop: semistable Cartier descent}, the isomorphism
\[\mathfrak{X}_{\eva}^{\assb}\cong \hig^{\assa}(C')
\]
in \Cref{prop: refinement of cartier descent} restricts to an isomorphism
\[
\mathfrak{X}_{\eva}^{\assb, \mathrm{ss}}\cong \hig^{\assa, \mathrm{ss}}(C').
\]
Combining with the isomorphism (\ref{eqn: restricted Li-Sun semistable}) and the equalities of parabolic degrees (\ref{eqn: matching parabolic degree}) and (\ref{eqn: matching degree of connections}), we obtain the desired isomorphism in \Cref{thm: main para}.
\end{proof}

\section{Cartier descent for filtered sheaves}\label{sec: filtered sheaf Cartier descent}
The goal of this section is to give a more explicit description of the isomorphism (\ref{eqn: equivalence between connection and parabolic bundle})
\[\con_0^{\asstau}(C)\cong \Bun^{\evalpha}(C', \evm')\]
in Proposition~\ref{prop: logarithmic Cartier descent}, showing that it may be regarded as a generalization of Cartier descent to filtered sheaves. We start by recalling the correspondence between parabolic bundles and locally free filtered sheaves as described in \cite{BodenYokogawa1996}. 
\begin{defn}
    A locally free filtered sheaf $\sE$ on $(C, q)$ consists of a locally free sheaf $\sE_0$ on $C$, together with an $\mathbb{R}$-indexed collection of locally free subsheaves $\{\sE_x\}_{x\in\mathbb{R}}\subset \sE_0(*q)$ satisfying the following four conditions:
 \begin{enumerate}
     \item\label{cond: inclusion} If $x\leq y$, then $\sE_{y}\subseteq \sE_{x}$,
      \item\label{cond: isomorphism away from q} For $x\leq y$, the inclusion $\sE_{y}\subseteq \sE_{x}$ is an
    isomorphism over $C\setminus q $,
     \item \label{cond: shift by 1}$\sE_{x+1}=\sE_x(-q)$,
     \item \label{cond:semicontinuous}
     For sufficiently small $\epsilon>0$, we have $\sE_{x-\epsilon}=\sE_x$.
 \end{enumerate}
 \end{defn}
 Given a parabolic bundle $(E, E_{\bullet}, \evalpha)$ with filtration\[0=E_{0}\subset E_{1}\subset\cdots\subset E_{r-1}\subset E_{r}=E_{q}\] and parabolic weights
 \[
\evalpha=(\alpha_1,\alpha_2,\ldots,\alpha_r), \qquad
1>\alpha_1>\alpha_2>\cdots>\alpha_r\geq 0,
\]
we set $\alpha_0=1$ and define
\[    \sE_{x}=\begin{cases}
E, & 0\leq x\leq \alpha_r,\\[2mm]
\ker(E\to E_q/E_j), & \alpha_{j+1}<x\leq \alpha_j,\qquad 0\leq j\leq r-1,
\end{cases}
\]
which extends to a unique collection $\{\sE_{x}\}_{x\in \mathbb{R}}$ satisfying the four conditions above. Conversely, given a locally free filtered sheaf $\{\sE_{x}\}_{x\in \mathbb{R}}$, we associate with it a parabolic bundle $(E, E_{\bullet}, \evalpha)$ as follows. We define $\on{gr}_x(\sE)\coloneqq\sE_x/\sE_{x+\epsilon}$. We set $E=\sE_0$, set
\[
\evalpha=\{x\in [0,1)\mid \on{gr}_x(\sE)\neq0\}=\{\alpha_1, \alpha_2,\dots,\alpha_r\}, \quad 1>\alpha_1>\alpha_2>\cdots>\alpha_r\geq 0,
\]
and set 
\[
E_j=\sE_{\alpha_j}/\sE_{1}\subseteq \sE_0/\sE_0(-q)\cong E_q \quad \text{for } j=1,2,\dots,r.
\]
It is straightforward to check that these constructions are mutually quasi-inverse functors, thus induce an equivalence between the groupoid of locally free filtered sheaves and the groupoid of parabolic bundles. 

\begin{defn}
    A filtered logarithmic connection $(\{\sE_x\}_{x\in \mathbb{R}}, \nabla)$ on $(C,q)$ consists of a locally free filtered sheaf $\{\sE_x\}_{x\in \mathbb{R}}$ together with a logarithmic connection 
    \[
    \nabla: \sE_0\longrightarrow \sE_0\otimes \omega_C(q)
    \]
    such that $\nabla(\sE_x)\subseteq \sE_x\otimes \omega_C(q)$ for each $x\in \mathbb{R}$, where $\nabla$ also denotes its natural extension to $\sE_0(*q)$.
\end{defn}
It is straightforward to check that the equivalence between locally free filtered sheaves and parabolic bundles described above upgrades to an equivalence between filtered logarithmic connections and parabolic connections.  

Denote by $\con_0^{\on{par}}(C)$ the moduli stack of parabolic connections on $(C,q)$ with zero $p$-curvature, and by $\Bun^{\on{par}}(C')$ the moduli stack of parabolic bundles on $(C',q')$. More precisely, we set
\[
\con_0^{\on{par}}(C)\coloneqq\coprod_{\evtau,\evbeta,\asstau}\con_0^{\asstau}(C),
\]
where the disjoint union ranges over all $\evtau\in \mathfrak{t}(\bF_p)/W$, parabolic weights $\evbeta$ and assignments $\asstau$; and we set
\[
\Bun^{\on{par}}(C')\coloneqq\coprod_{\evalpha,\evm'}\Bun^{\evalpha}(C',\evm'),
\]
where the disjoint union ranges over all ordered tuples $\evm'$ that sum to $n$, and sets of parabolic weights $\evalpha$ such that $\#\evalpha$ is equal to the length of $\evm'$.

Let $\cE=(E,E_{\bullet}, \nabla_E)\in \con_0^{\on{par}}(C)$, and let $(\{\sE_x\}_{x\in \mathbb{R}},\nabla_{\sE})$ be the filtered connection corresponding to $\cE$. For any $x\in \mathbb{R}$, we define
\begin{equation}\label{eqn: filtered sheaf cartier descent}
\sF_{x}\coloneqq(\on{Fr}_C)_*\ker\bigl(
\nabla_{\sE}\colon \sE_{px}
\longrightarrow
\sE_{px}\otimes\omega_C(q)
\bigr).
\end{equation}
\begin{lemma}\label{lemma: kernel gives filtered sheaf}
    The $\mathbb{R}$-collection $\{\sF_{x}\}_{x\in\mathbb{R}}$ is a locally free filtered sheaf over $(C',q')$. 
\end{lemma}
\begin{proof}
    We start by showing $\sF_x$ is locally free for each $x$. Recall that by \Cref{lemma: local isomorphism}, $\cE$ is locally isomorphic to the model object $\cE^0$ defined in \Cref{sec: proof of identifying with parabolic higgs}. The fact that each $\sF_x$ is locally free follows from the simple computation that 
\[
(\Fr_C)_*\ker\bigl(
\nabla^0\colon \cO_C(\gamma q)
\longrightarrow
\cO_C(\gamma q)\otimes\omega_C(q)
\bigr)
=
\cO_{C'}\bigl(\lfloor\gamma/p\rfloor q'\bigr)
\] 
for any integer $\gamma$. This argument works for $\cE$ over an arbitrary test $k$-scheme. 

The conditions (\ref{cond: inclusion}), (\ref{cond: isomorphism away from q}), (\ref{cond:semicontinuous}) follow from the fact that $\{\sE_x\}_{x\in \mathbb{R}}$ satisfies those conditions. Let $x$ be a local coordinate at $q$. Condition (\ref{cond: shift by 1}) follows from the observation that for any section $s$ of $\sE_0(*q)$ around $q$, we have
    \[    \nabla_{\sE}(s)=0\Longleftrightarrow \nabla_{\sE}(x^ps)=0.
    \]
\end{proof}
By Lemma~\ref{lemma: kernel gives filtered sheaf}, the $\mathbb{R}$-collection of sheaves
$\{\sF_x\}_{x\in\mathbb{R}}$ defined in
\eqref{eqn: filtered sheaf cartier descent} forms a locally free
filtered sheaf on $(C',q')$. Via the correspondence with parabolic
bundles, this construction defines a morphism of stacks
\[
\Phi:\con_0^{\on{par}}(C)\longrightarrow\Bun^{\on{par}}(C').
\]

\begin{thm}\label{thm: filtered sheaf cartier descent}
    For any $\evtau\in \mathfrak{t}(\bF_p)/W$, parabolic weights $\evbeta$ and assignment $\asstau$, the morphism $\Phi$ restricts to an isomorphism
    \[\con_0^{\asstau}(C)\cong \Bun^{\evalpha}(C', \evm')\]
that agrees with the one constructed in Proposition~\ref{prop: logarithmic Cartier descent}. Since \eqref{eqn: asstau to assa} gives a bijection between assignments $\asstau$ and pairs $(\evalpha,\evm')$, $\Phi$ induces an isomorphism
\[
\Phi: \con_0^{\on{par}}(C)\xrightarrow{\sim}\Bun^{\on{par}}(C').
\]
\end{thm}
\begin{proof}
    We denote the isomorphism constructed in Proposition~\ref{prop: logarithmic Cartier descent} by
    \[\Phi'_{\asstau}:\con_0^{\asstau}(C)\cong \Bun^{\evalpha}(C', \evm')\] 
and the restriction of $\Phi$ to $\con_0^{\asstau}(C)$ by $\Phi_{\asstau}$. Let $\cE^0$ be the model object defined in Section \ref{sec: proof of identifying with parabolic higgs}. We have the following $2$-commutative diagram
\[
\begin{tikzcd}[column sep=huge]
\Bun_{{\cA}ut(\cE^0)}(C)\arrow[r, "\text{Lemma}~\ref{lemma: connection auto and parabolic group scheme}", "\sim"'] \arrow[d, "\sim","\mathsf{G}"'] & \Bun_{\cG_{\widetilde{\evalpha}}}(C') \arrow[d, "\text{Remark}~\ref{remark: parabolic bundle and parabolic group torsor}", "\sim"'] \\
\con_0^{\asstau}(C) \arrow[r, "\sim","\Phi'_{\asstau}"'] & \Bun^{\evalpha}(C', \evm'),
\end{tikzcd}
\]
where $\mathsf{G}$ is constructed in the proof of Proposition~\ref{prop: logarithmic Cartier descent}. Recall that for any $\cP\in \Bun_{{\cA}ut(\cE^0)}(C)$, we have
\[
\mathsf{G}(\cP)=\cP\times^{{\cA}ut(\cE^0)}\cE^0,
\]
with connection and parabolic structure inherited from $\cE^0$. Let $(\{\sE_x^0\}_{x\in \bR}, \nabla^0_{\sE})$ be the filtered connection corresponding to $\cE^0$. Then the filtered sheaf corresponding to $\mathsf{G}(\cP)$ is given by \[\cP\times^{{\cA}ut\bigl((\{\sE_x^0\}_{x\in \bR}, \nabla^0_{\sE})\bigr)}\{\sE_x^0\}_{x\in \bR}. 
\]
By the definition of $\Phi$ through (\ref{eqn: filtered sheaf cartier descent}), we have 
\[
\Phi_{\asstau}(\mathsf{G}(\cP))\cong(\on{Fr}_C)_*\cP\times^{(\on{Fr}_C)_*{\cA}ut(\cE^0)}\Phi(\cE^0).
\]
By the $2$-commutative diagram above, we have 
\[
\Phi'_{\asstau}(\mathsf{G}(\cP))\cong (\on{Fr}_C)_*\cP\times^{\cG_{\widetilde{\evalpha}}}\cF^0,
\]
where $\cF^0=(F^0,F^0_{\bullet})$ is the model object constructed in Remark~\ref{remark: parabolic bundle and parabolic group torsor}. Therefore in order to show $\Phi_{\asstau}\simeq \Phi'_{\asstau}$, it is enough to identify $\Phi(\cE^0)\cong \cF^0$, such that under this identification, the action of $(\on{Fr}_C)_*{\cA}ut(\cE^0)$ on $\Phi(\cE^0)$ agrees with the action of $\cG_{\widetilde{\evalpha}}$ on $\cF^0$ through the canonical isomorphism
\[
(\on{Fr}_C)_*{\cA}ut(\cE^0)\cong \cG_{\widetilde{\evalpha}}
\]
in Lemma~\ref{lemma: connection auto and parabolic group scheme}. Let $\cE^0_i=\bigl(\cO_{C}(-\nu(\tau_i)q), \nabla^0,\beta_{\on{op}}^{\widetilde{\evtau}}(i)\bigr)$ be the rank one parabolic connection with parabolic weight $\beta_{\on{op}}^{\widetilde{\evtau}}(i)$, with the logarithmic connection $\nabla^0$ induced from the canonical connection on $\cO_C$. Then we have $\cE^0=\bigoplus_{i=1}^n \cE^0_i$. Note that $\Phi$ respects direct sums, therefore the desired identification $\Phi(\cE^0)\cong \cF^0$ follows from the rank one computation that 
\[
\Phi(\cE^0_i)\cong(\cO_{C'}, \alphaop(i)), \quad \text{where } \alphaop(i)=\frac{\beta_{\on{op}}^{\widetilde{\evtau}}(i)-\nu(\tau_i)}{p},
\]
see Lemma~\ref{lemma: rank 1 cartier descent}.
The agreement of actions follows from the computation in Lemma~\ref{lemma: connection auto and parabolic group scheme}.
\end{proof}

\begin{lemma}\label{lemma: rank 1 cartier descent}
    For any $\beta\in[0,1)$ and $\gamma\in\{0,1,2,\dots,p-1\}$, let $\alpha=\frac{\beta+\gamma}{p}$, then we have
    \[
    \Phi\bigl((\cO_C(\gamma q), \nabla^0, \beta)\bigr)\cong (\cO_{C'}, \alpha).
    \]
\end{lemma}
\begin{proof}
    Let $(\{\sE_x\}_{x\in\mathbb{R}},\nabla_{\sE}^0)$ be the filtered connection corresponding to $(\cO_C(\gamma q), \nabla^0, \beta)$ and let $\{\sF_x\}_{x\in\mathbb{R}}$ be the filtered sheaf corresponding to $\Phi\bigl((\cO_C(\gamma q), \nabla^0, \beta)\bigr)$. By definition, we have
    \[    \sF_{x}=(\on{Fr}_C)_*\ker\left(
\nabla^0_{\sE}\colon \sE_{px}
\longrightarrow
\sE_{px}\otimes\omega_C(q)
\right)
    \]
    as in (\ref{eqn: filtered sheaf cartier descent}). Setting $x=0$, we have 
    \[   \sF_{0}=(\on{Fr}_C)_*\ker\left(
\nabla^0\colon \cO_C(\gamma q)
\longrightarrow
\cO_C(\gamma q)\otimes\omega_C(q)
\right)=\cO_{C'}.
    \]
Note that the only $x\in[0,1)$ where $\sF$ jumps is the one where $\sE_{px}$ jumps from $\cO_C$ to $\cO_C(-q)$, which happens when $px=\beta+\gamma$.
\end{proof}

\begin{remark}\label{Remark: inverse of Phi}
The inverse of 
\[
\Phi: \con_0^{\on{par}}(C)\xrightarrow{\sim}\Bun^{\on{par}}(C')
\]
can be described as follows. Let $\{\sF{_x}\}_{x\in{\mathbb{R}}}$ be a locally free filtered sheaf on $(C',q')$. We define
\[
\sE_{x}\coloneqq \sum_{\gamma=0}^{p-1} (\on{Fr}_C^*\sF_{(x+\gamma)/p})(\gamma q)\subset (\on{Fr}_C^*\sF_0)(*q).
\]
 It is clear that $\{\sE_x\}_{x\in \mathbb{R}}$ is a locally free filtered sheaf on $(C,q)$. The restriction of the canonical connection on $(\on{Fr}_C^*\sF_0)(*q)$ gives rise to a connection $\nabla^0$ on each $\sE_x$, and we have $(\{\sE_x\}_{x\in \mathbb{R}}, \nabla^0)\in \con_0^{\on{par}}(C)$. Let 
\[
\Psi: \Bun^{\on{par}}(C')\longrightarrow \con_0^{\on{par}}(C)
\]
be the morphism of stacks mapping $\{\sF{_x}\}_{x\in\mathbb{R}}$ to $(\{\sE_x\}_{x\in \mathbb{R}}, \nabla^0)$. Note that 
\[
\on{Fr}_C^*\sF_{x}\subseteq\sum_{\gamma=0}^{p-1} (\on{Fr}_C^*\sF_{(px+\gamma)/p})(\gamma q)\subseteq (\on{Fr}_C^*\sF_{x})((p-1)q),
\]
therefore
\[
(\on{Fr}_C)_*\ker\left(\nabla^0|_{\on{Fr}_C^*\sF_{x}}\right)\subseteq \left(\Phi\circ \Psi (\{\sF_x\}_{x\in\mathbb{R}})\right)_x\subseteq (\on{Fr}_C)_*\ker\left(\nabla^0|_{(\on{Fr}_C^*\sF_{x})((p-1)q)}\right).
\]
We have 
\[
(\on{Fr}_C)_*\ker\left(\nabla^0|_{\on{Fr}_C^*\sF_{x}}\right)=(\on{Fr}_C)_*\ker\left(\nabla^0|_{(\on{Fr}_C^*\sF_{x})((p-1)q)}\right)=\sF_x,
\]
therefore 
\[
\left(\Phi\circ \Psi (\{\sF{_x}\}_{x\in\mathbb{R}})\right)_x=\sF_x.
\]
\end{remark}
\begin{remark}
    The formulation of Cartier descent for filtered sheaves through (\ref{eqn: filtered sheaf cartier descent}) potentially extends to higher dimensional varieties. We plan to investigate this generalization in future work.
\end{remark}
\newpage

\bibliographystyle{alpha}
\footnotesize{\bibliography{parabolic.bib}}
\end{document}